\documentclass[a4paper]{article}

\pdfoutput=1

\usepackage[utf8]{inputenc}
\usepackage[T1]{fontenc}
\usepackage[UKenglish]{babel}
\usepackage{microtype}
\usepackage{amsmath,amssymb,amsthm}
\usepackage{hyperref}
\usepackage{cleveref}
\usepackage{booktabs}
\usepackage{makecell}
\usepackage{tikz}
\usepackage{graphicx}
\usepackage[margin=1in]{geometry}
\usepackage[font=small,skip=8pt]{caption}

\usetikzlibrary{shapes.geometric, arrows.meta, positioning, calc}

\newtheorem{theorem}{Theorem}

\newtheorem{remark}[theorem]{Remark}

\title{Adaptive Subproblem Selection in Benders Decomposition\\for Survivable Network Design Problems}

\author{Tim Donkiewicz\\
Chair of Operations Research, RWTH Aachen University, Germany\\
\texttt{donkiewicz@or.rwth-aachen.de}}

\date{}

\newcommand{\R}{\mathbb{R}}
\newcommand{\Z}{\mathbb{Z}}
\newcommand{\N}{\mathbb{N}}
\newcommand{\todo}[1]{}
\newcommand{\fixme}[1]{}
\newcommand{\question}[1]{}
\newcommand{\critical}[1]{}
\newcommand{\note}[1]{}

\begin{document}

\maketitle

\begin{abstract}
Scenario-based optimization problems can be solved via Benders decomposition, which separates first-stage (master problem) decisions from second-stage (subproblem) recourse actions and iteratively refines the master problem with Benders cuts.
In conventional Benders decomposition, all subproblems are solved at each iteration.
For problems with many scenarios, solving only a selected subset can reduce computation.
We quantify the potential in selecting only those subproblems that yield cuts, and develop subproblem scoring and selection strategies. The proposed multi-criteria scoring methods combine historical subproblem performance metrics with problem-specific features, trained online via logistic regression to adapt to the changing likelihood of subproblem usefulness. Multiple stopping criteria balance exploration and exploitation: cut limits, proportional solve limits, and score thresholds. We evaluate our approach on a variant of the survivable network design problem, which serves as a testbed due to its natural decomposition into many subproblems of varying importance.

Computational experiments on 135 test instances demonstrate the potential and practical performance of subproblem selection. Analysis reveals that 52.1\% of all subproblems solved are unnecessary (they contribute no cuts and occur outside cut-free rounds).
An oracle with perfect foresight reduces total solve times by 34.4\%. Random selection performs significantly worse than full enumeration, showing that naive strategies can degrade performance. Our best-scoring and selection method achieves statistically significant improvements in both runtime and primal-dual integrals. These results provide empirical evidence that informed subproblem selection can improve Benders decomposition in this setting, while highlighting challenges in developing reliable prediction models. Whether these findings extend to other problem classes remains an open question for future work.
\end{abstract}


\section{Introduction}
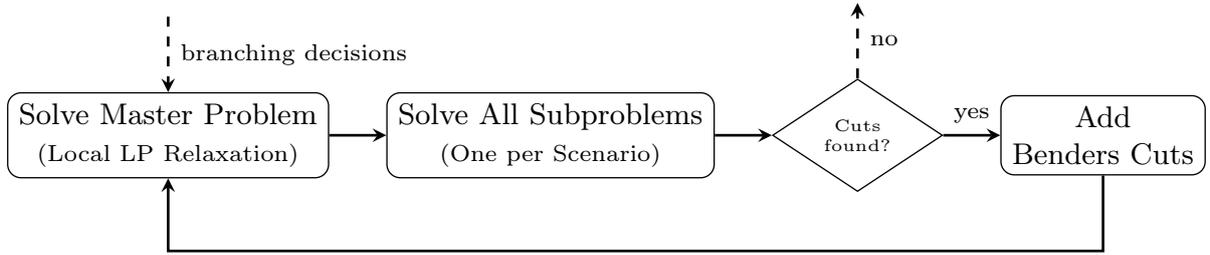
\begin{figure}[htbp]
    \centering
    \resizebox{\textwidth}{!}{\begin{tikzpicture}[
    node distance=0.3cm and 0.6cm,
    box/.style={rectangle, draw, rounded corners, minimum width=1.8cm, minimum height=0.8cm, align=center, font=\small},
    decision/.style={diamond, draw, minimum width=1.3cm, minimum height=0.55cm, align=center, font=\tiny, aspect=1.5},
    arrow/.style={->, >=stealth, thick}
]

\node[box] (master) {Solve Master Problem\\\scriptsize(Local LP Relaxation)};
\node[box, right=of master] (subproblems) {Solve All Subproblems\\\scriptsize(One per Scenario)};
\node[decision, right=of subproblems] (cutsfound) {Cuts\\found?};
\node[box, right=of cutsfound] (cuts) {Add\\Benders Cuts};

\draw[arrow] (master) -- (subproblems);
\draw[arrow] (subproblems) -- (cutsfound);
\draw[arrow] (cutsfound) -- node[above, font=\scriptsize] {yes} (cuts);

\draw[arrow, dashed] ([yshift=0.8cm]master.north) -- node[right, font=\scriptsize] {branching decisions} (master.north);

\draw[arrow] (cuts.south) -- ++(0,-0.8) -| (master.south);

\draw[arrow, dashed] (cutsfound.north) -- node[right, font=\scriptsize] {no} ++(0,0.8);

\end{tikzpicture}}
    \caption{Basic Branch-and-Benders-cut framework.}
    \label{fig:branch-and-benders-cut}
\end{figure}

Benders decomposition \cite{benders1962partitioning} is a powerful technique for structured optimization problems, separating first-stage decisions from second-stage recourse problems. In an iterative manner, the master problem is solved to propose first-stage solutions. Then, subproblems corresponding to different scenarios are solved to check feasibility or optimality of the proposed solution. If any subproblem is infeasible or more expensive than expected by the master problem's estimation, Benders cuts are generated and added to the master problem to refine the solution space. This process repeats until the branch-and-bound algorithm terminates.
Modern implementations of Benders decomposition apply Branch-and-Benders-cut \cite{rahmaniani2017benders}, integrating Benders cuts as lazy constraints within a branch-and-bound framework for mixed-integer programming (see \autoref{fig:branch-and-benders-cut}).

However, when the number of scenarios grows large, solving all subproblems at each iteration becomes computationally prohibitive. Through computational experiments on survivable network design (SND) instances, we observe that 52.1\% of subproblems solved under conventional Benders decomposition are unnecessary. This empirical observation motivates a data-driven computational approach to adaptive subproblem selection---we aim to dynamically learn which subproblems are most likely to contribute to the solving process by generating Benders cuts. To that end, we observe patterns in cut generation behavior, then prioritize subproblems accordingly, and stop once enough subproblems were solved (and at least one cut was generated).

We use SND as a testbed for developing adaptive subproblem selection strategies. SND problems exhibit characteristics that make them particularly suitable for this study: they naturally decompose into many subproblems (one per failure scenario), these subproblems vary substantially in their importance, and some scenarios can dominate others due to network topology. While our goal is not to advance the state-of-the-art for SND solving specifically, this problem class provides a suitable environment to demonstrate subproblem selection techniques.

To our knowledge, adaptive subproblem selection based on computational pattern learning has not been successfully demonstrated for Benders decomposition in prior work. We provide a first empirical study addressing this gap, using SND as our testbed. Our four main contributions are grounded in empirical analysis: First, we develop an oracle baseline that quantifies the computational potential of perfect subproblem selection through empirical measurements. Second, we develop a multi-criteria scoring mechanism that combines observed historical performance metrics (cut generation frequency, cumulative contribution, recency) with computational features extracted from the current solution (failed edge capacity, flow, utilization, centrality). Third, we introduce an online learning approach using logistic regression that learns from computational observations during solving. We train incrementally on each subproblem outcome to predict scenario infeasibility and adapt to evolving cut generation patterns without requiring problem-specific tuning or offline training. The online training incurs a reasonable computational overhead---one update is comparable to solving a single subproblem---making it practical for real-time deployment. Fourth, we implement partial subproblem solving with multiple stopping criteria that balance computational effort between exploration and exploitation based on empirical performance indicators.

Our goal is twofold: demonstrate through computational experiments on SND instances that dynamically learning from empirical observations can achieve statistically significant improvements over the baseline, and quantify via the oracle how much additional computational improvement could be attained with more sophisticated pattern recognition. Validating these findings on other problem classes remains a direction for future research. 

\section{Related Work}
\label{sec:related-work}

Benders \cite{benders1962partitioning} first introduced the Benders decomposition algorithm for mixed-integer programming. Rahmaniani et al.~\cite{rahmaniani2017benders} provide a comprehensive modern survey and an overview of Branch-and-Benders-cut. Benders decomposition has been applied to survivable network design problems \cite{botton2013benders,fortz2009improved}. The SNDlib benchmark instances \cite{IdzikowskiOrlowskiRaackWoesnerWolisz2010b,KosterKutschkaRaack2010} are widely used for evaluating network design algorithms. 

\subsection{Scenario Retention, Selection, and Partitioning}
\label{subsec:scenario-selection}

An approach to managing computational burden is to retain selected second-stage variables explicitly in the master problem rather than projecting them out through Benders cuts. Crainic et al.~\cite{crainic2014partial,crainic2021partial} develop partial decomposition strategies for two-stage integer programs and stochastic multicommodity network design. Pauphilet et al.~\cite{pauphilet2025random} show that random variable retention can achieve good performance. These methods fundamentally change the master problem structure.

Rather than modifying the master problem structure, scenario selection methods choose which subproblems to solve at each iteration. Rahmaniani et al.~\cite{rahmaniani2017benders} note that solving only a subset of subproblems---especially early in the algorithm---can reduce computational burden, though such strategies remain underexplored for combinatorial optimization. Chen et al.~\cite{chen2011implicit} reformulate the feasibility check for single-commodity survivable network design using an oracle based on minimum cut computations to identify violated scenarios to examine. Zhang et al.~\cite{zhang2023optimized} propose optimization-based scenario reduction to approximate the recourse function while maintaining tight optimality gaps. Blanchot et al.~\cite{blanchot2023benders} stop solving subproblems after the sum of the violation of obtained optimality cuts exceeds the current gap, often solving only 1\% of subproblems. Mehamdi and Gendron~\cite{mehamdi2024pricing} select pricing subproblems to solve in partial pricing in Branch-and-Price based on dual information. Celik and Toriello~\cite{celik2025exact} develop a scenario-retention strategy that avoids resolving subproblems remaining feasible across iterations. Learning-based approaches have also emerged: Borozan et al.~\cite{borozan2024machine} use machine learning to predict which scenarios generate cuts binding in the final solution. These approaches adaptively select which subproblems to examine without modifying the problem structure.

\subsection{Scoring Mechanisms for Algorithmic Decisions}
\label{subsec:scoring-mechanisms}

The use of scoring functions to prioritize algorithmic decisions is common in mixed-integer programming \cite{achterberg2007constraint,achterberg2009scip,andreello2007embedding}, and multi-criteria scoring methods often outperform, e.g., pure violation-based approaches \cite{hosseini2024deepest,avella2004metric}.
Recent work has also explored machine learning for various selection decisions in MIP solvers, applying neural networks and reinforcement learning to cut selection, variable selection, and cut generation, learning to predict effectiveness from historical data and solution trajectories \cite{scavuzzo2024machine,tang2020reinforcement,turner2023adaptive,wang2023learning,zhang2025learning}. Li et al.~\cite{li2023learning} use machine learning to configure separator selection in branch-and-cut. Deza et al.~\cite{deza2023machine} provide a comprehensive survey showing that while machine learning approaches demonstrate promise, they require substantial training data and careful feature engineering.

These principles from MIP solver design inform our multi-criteria scoring for subproblem selection, where we decide which feasibility checks to perform rather than selecting cuts from an existing pool.

\subsection{Contributions in the Context of Related Works}
\label{subsec:positioning}

Prior work on subproblem selection designs auxiliary optimization problems or heuristics to approximate subproblem objectives and identify scenarios to examine, often leveraging problem-specific structure. Our approach differs by scoring subproblems using multiple criteria and adaptively determining how many to solve without requiring problem reformulation. Through online learning from empirical observations during solving, we dynamically adapt prioritization to each instance and solution trajectory. We observe that for our test instances with feasibility cuts only, omitting even few cuts significantly degrades convergence, highlighting the challenge of effective subproblem selection.

Our multi-criteria scoring framework adapts principles from MIP solver decision-making (\autoref{subsec:scoring-mechanisms}) to decide which subproblems to solve before incurring costs, reducing computational effort by avoiding unnecessary solves entirely.

We do not aim to advance the state-of-the-art for survivable network design specifically, but rather use it as a representative problem class to develop and evaluate adaptive subproblem selection strategies applicable to Benders decomposition more broadly.

\section{Problem Formulation}
\label{sec:problem}

Consider a network $G = (V, E)$ with node set $V$ and undirected edge set $E$.
For each edge $\{i,j\} \in E$, we create two directed arcs: $(i,j)$ and $(j,i)$.
Let $A$ denote the set of all such directed arcs. For any node $v \in V$, we denote by $\delta^+(v)$ the set of arcs leaving $v$ and by $\delta^-(v)$ the set of arcs entering $v$.
We denote by $A(e)$ the two arcs corresponding to edge $e \in E$.

For each edge $e \in E$, capacity can be installed using a set of \emph{capacity modules} $M_e$.
Each module type $m \in M_e$ provides capacity $u_{e,m} > 0$ at installation cost $c_{e,m} \geq 0$.
Edges may have preinstalled capacity (at zero cost).
The number of installed modules $m$ on each edge $e$ is represented by $y_{e,m} \in \N$.

A set of demands $D$ must be routed, where each demand $d \in D$ specifies a source node $o_d\in V$, destination node $t_d \in V$, and demand value $b_d > 0$. For convenience, define $b_{d,v} = b_d$ if $v = o_d$, $b_{d,v} = -b_d$ if $v = t_d$, and $b_{d,v} = 0$ otherwise.

We consider a set $S$ of \emph{failure scenarios}, where each scenario $s \in S$ is characterized by a subset $F_s \subseteq E$ of failed edges. We also define $s_0$ as the \emph{base scenario} with $F_{s_0} = \emptyset$ (no failures), and let $S_0 = \{s_0\} \cup S$ denote all scenarios including the base case. The formulation supports arbitrary failure sets (single-edge, multi-edge, or scenario-specific combinations). In this paper, we focus on \emph{single-edge failure scenarios}: $S = \{s_e : e \in E\}$ where $F_{s_e} = \{e\}$ (edge $e$ fails), commonly known as N-1 reliability in telecommunications and power systems. The generalization to arbitrary failure sets is straightforward and does not affect the adaptive subproblem selection methodology.

This problem can be viewed as a two-stage integer program where, in each failure scenario, the demands must be routed.
In the {first stage}, we install capacity modules on edges (integer decisions made before failures occur) and determine flow routing for the base scenario $s_0$.
Then, in the {second stage}, for each scenario $s \in S$, we route demands on non-failed edges respecting installed capacities (recourse decisions after failures are revealed).
The {objective} is to minimize total module installation costs.

Note that the second-stage decisions (flow routing) serve as recourse actions but have zero cost---we only consider routing feasibility, not routing costs. The objective is purely to minimize first-stage installation costs, but the installed capacity must be sufficient to enable feasible routing under all failure scenarios.

\subsection{Compact Formulation}

Let $y_{e,m} \in \Z_{\geq 0}$ be the number of modules of type $m$ installed on link $e$ and $f_{s,d,a} \in \R_+$ the flow of demand $d$ on arc $a$ in scenario $s$.
We then consider the following compact formulation:
{\footnotesize
\begin{align}
    \min_{y} \quad & \sum_{e \in E} \sum_{m \in M_e} c_{e,m} y_{e,m} \label{eq:compact_obj} \\
    \text{s.t.} \quad & \sum_{a \in \delta^+(v)} f_{s,d,a} - \sum_{a \in \delta^-(v)} f_{s,d,a} = b_{d,v} && \forall s \in S_0, d \in D, v \in V \label{eq:flow_cons} \\
    & \sum_{d \in D} \sum_{a \in A(e)} f_{s,d,a} \leq \sum_{m \in M_e} u_{e,m} y_{e,m} && \forall s \in S_0, e \in E \setminus F_s \label{eq:capacity} \\
    & f_{s,d,a} = 0 && \forall s \in S_0, d \in D, e \in F_s, a \in A(e) \label{eq:failed_flow} \\
    & y_{e,m} \in \Z_{\geq 0} && \forall e \in E, m \in M_e \label{eq:vars_y} \\
    & f_{s,d,a} \geq 0 && \forall s \in S_0, d \in D, a \in A \label{eq:vars_f}
\end{align}
}

The capacity constraints~\eqref{eq:capacity} apply only to operational links $e \in E \setminus F_s$ in each scenario. For failed links $e \in F_s$, constraints~\eqref{eq:failed_flow} explicitly force flow to zero on both arcs $a \in A(e)$ of the failed link, prohibiting any routing through failed components.

With $|S_0| = |E| + 1$ scenarios (base case plus all single-edge failures), we have $O(|S_0| \cdot |D| \cdot |A|)$ flow variables and $O(|S_0| \cdot (|D| \cdot |V| + |E|))$ constraints. When considering multi-edge failures, the number of scenarios grows exponentially, making the compact formulation impractical. For large networks, even solving the single-edge failure case becomes computationally intractable using conventional MIP solvers.

\subsection{Benders Decomposition}

Benders decomposition addresses the scalability issue by separating first-stage and second-stage decisions into a master problem and subproblems.
Since we only consider feasibility in the second stage, we only use Benders \emph{feasibility} cuts to iteratively refine the master problem.

\paragraph*{Master Problem}

The Benders master problem includes first-stage capacity decisions and base case flow constraints:
{\footnotesize
\begin{align}
    \min_{y} \quad & \sum_{e \in E} \sum_{m \in M_e} c_{e,m} y_{e,m} \label{eq:master_obj} \\
    \text{s.t.} \quad & \sum_{a \in \delta^+(v)} f_{s_0,d,a} - \sum_{a \in \delta^-(v)} f_{s_0,d,a} = b_{d,v} && \forall d \in D, v \in V \label{eq:master_flow} \\
    & \sum_{d \in D} \sum_{a \in A(e)} f_{s_0,d,a} \leq \sum_{m \in M_e} u_{e,m} y_{e,m} && \forall e \in E \label{eq:master_capacity} \\
    & \text{Benders cuts \eqref{eq:benders_cut}} \label{eq:benders_cuts} \\
    & y_{e,m} \in \Z_{\geq 0} && \forall e \in E, m \in M_e \label{eq:master_vars_y} \\
    & f_{s_0,d,a} \geq 0 && \forall d \in D, a \in A \label{eq:master_vars_f}
\end{align}
}
Constraints \eqref{eq:master_flow} and \eqref{eq:master_capacity} enforce flow conservation and capacity constraints for the base case scenario $s_0$ with no failed edges ($F_{s_0} = \emptyset$), ensuring basic feasibility before considering failure scenarios.

\paragraph*{Subproblem}

For each failure scenario $s \in S$ and candidate solution $\bar{y}$, the subproblem checks feasibility:
{\footnotesize
\begin{align}
    \min \quad & 0 \label{eq:sub_obj} \\
    \text{s.t.} \quad & \sum_{a \in \delta^+(v)} f_{s,d,a} - \sum_{a \in \delta^-(v)} f_{s,d,a} = b_{d,v} && \forall d \in D, v \in V \label{eq:sub_flow} \\
    & \sum_{d \in D} \sum_{a \in A(e)} f_{s,d,a} \leq \sum_{m \in M_e} u_{e,m} \bar{y}_{e,m} && \forall e \in E \setminus F_s \label{eq:sub_capacity} \\
    & f_{s,d,a} = 0 && \forall d \in D, e \in F_s, a \in A(e) \label{eq:sub_failed_flow}
\end{align}
}

Constraints \eqref{eq:sub_flow} enforce flow conservation (corresponding to \eqref{eq:flow_cons} in the compact formulation), and constraints \eqref{eq:sub_capacity} enforce capacity limits on operational edges (corresponding to \eqref{eq:capacity}). For failed edges $e \in F_s$, constraints \eqref{eq:sub_failed_flow} explicitly set flow to zero on both arcs, ensuring no routing through failed components.

We observe that this is a multicommodity flow problem with fixed capacities determined by the first-stage solution $\bar{y}$.
It is solvable in polynomial time as a linear program with $O(|D| \cdot |A|)$ variables and $O(|D| \cdot |V| + |E|)$ constraints.

If the subproblem is feasible, the capacity $\bar{y}$ is sufficient for scenario $s$.
If infeasible, we generate a Benders feasibility cut using the subproblem's dual solution.

\paragraph*{Benders Feasibility Cuts}

When the subproblem for scenario $s$ is infeasible, we obtain an unbounded ray $(\pi^s, \sigma^s)$ of the dual problem, where $\pi_{e}^s \geq 0$ are the dual variables for capacity constraints \eqref{eq:sub_capacity} and $\sigma_{d,v}^s$ are the (unrestricted) dual variables for flow conservation constraints \eqref{eq:sub_flow}. The Benders feasibility cut is:
\begin{equation}
    \sum_{e \notin F_s} \pi_{e}^s \sum_{m \in M_e} u_{e,m} y_{e,m} \geq -\sum_{d \in D} \sum_{v \in V} \sigma_{d,v}^s b_{d,v}
    \label{eq:benders_cut}
\end{equation}
where we sum only over operational edges $e \notin F_s$ since failed edges have zero available capacity.

\section{Adaptive Subproblem Selection}
\label{sec:adaptive-filtering}

Our data-driven approach is based on two observations.
Firstly, it is not required to solve all subproblems at each Benders iteration to make progress toward feasibility.
Secondly, not all subproblems are equally likely to generate cuts.
If we are able to predict which scenarios are unlikely to produce a cut, we can avoid solving those subproblems entirely, saving computation time.

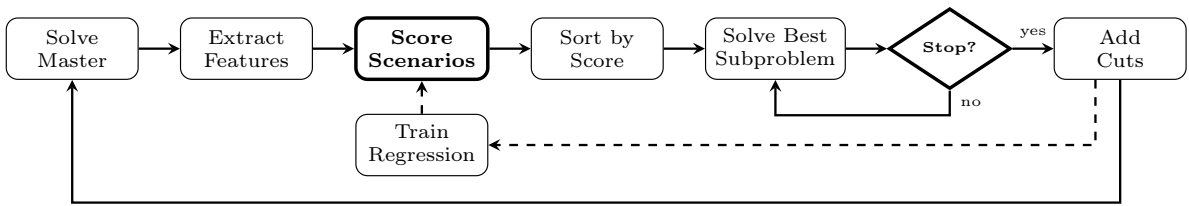
\begin{figure}[ht]
\centering
\resizebox{\textwidth}{!}{
\begin{tikzpicture}[
    node distance=0.3cm and 0.5cm,
    box/.style={rectangle, draw, rounded corners, minimum width=1.6cm, minimum height=0.75cm, align=center, font=\scriptsize},
    boxbold/.style={rectangle, draw, line width=1.2pt, rounded corners, minimum width=1.6cm, minimum height=0.75cm, align=center, font=\scriptsize\bfseries},
    decision/.style={diamond, draw, minimum width=1.3cm, minimum height=0.55cm, align=center, font=\tiny, aspect=1.5},
    decisionbold/.style={diamond, draw, line width=1.2pt, minimum width=1.3cm, minimum height=0.55cm, align=center, font=\tiny\bfseries, aspect=1.5},
    arrow/.style={->, >=stealth, thick}
]

\node[box] (master) {Solve\\Master};
\node[box, right=of master] (extract) {Extract\\Features};
\node[boxbold, right=of extract] (score) {Score\\Scenarios};
\node[box, right=of score] (sort) {Sort by\\Score};
\node[box, right=of sort] (solve) {Solve Best\\Subproblem};
\node[decisionbold, right=of solve] (stop) {Stop?};
\node[box, right=of stop] (add) {Add\\Cuts};

\node[box, below=0.4 of score] (train) {Train\\Regression};

\draw[arrow] (master) -- (extract);
\draw[arrow] (extract) -- (score);
\draw[arrow] (score) -- (sort);
\draw[arrow] (sort) -- (solve);
\draw[arrow] (solve) -- (stop);
\draw[arrow] (stop) -- node[above, font=\tiny] {yes} (add);

\draw[arrow] (stop.south) -- node[right, font=\tiny] {no} ++(0,-0.3) -| (solve.south);

\draw[arrow] (add.south) -- ++(0,-1.5) -| (master.south);

\draw[arrow, dashed] ([xshift=-0.3cm]add.south) -- ++(0,-0.5) |- (train.east);
\draw[arrow, dashed] (train) -- (score);

\end{tikzpicture}
}
\caption{Adaptive subproblem selection pipeline in Benders decomposition. Scoring and subproblem selection step are highlighted.}
\label{fig:pipeline}

\end{figure}

We therefore propose Branch-and-Benders-cut with \emph{adaptive subproblem selection}: prioritize which subproblems to solve first based on learned importance, and limit the number of subproblems solved per iteration, as illustrated in \autoref{fig:pipeline}.
At each Benders iteration, we first extract features from the master solution (topology metrics like capacity, flow, and utilization).
We then score all scenarios based on their predicted likelihood of generating cuts using a machine learning model.
Scenarios are examined in decreasing score order, with stopping criteria determining when to terminate subproblem solving. If there exists any, at least one cut is added in each round. After each subproblem is solved, we train the regression model online using the observed outcome and updated historical performance (reliability, violation magnitude). Feasibility cuts from infeasible subproblems are added to the master, and the process repeats until convergence.

In this section, we first describe an oracle that selects an optimal subproblem solving sequence to not miss any cuts.
Then, we present our multi-criteria subproblem scoring mechanism in Section~\ref{sec:subproblem-scoring}.
Building upon the scoring system, we introduce stopping criteria for subproblem solving iterations in Section~\ref{sec:partial-subproblem-solving}, aiming to approximate the oracle's behavior.

\subsection{Oracle for Evaluating Prediction Quality}

To assess scoring quality, we define a \emph{cut prediction oracle} $\mathcal{O}^*$ with perfect foresight of all cuts generated until the algorithm converges. Let $\Pi = (P_1, P_2, \ldots, P_T)$ denote a subproblem selection sequence, where $P_t \subseteq S$ are the examined scenarios in iteration $t$. We define the set of cut-yielding subproblems in iteration $t$ as $V_t \subseteq S$.
The \textit{optimal subproblem selection} as given by $\mathcal{O}^*$ is $\Pi^*=(V_1, V_2, \ldots, V_t)$, the sequence consisting of exactly those subproblems yielding a cut.
Note that a sequence $\Pi$ with $P_t \cap V_t \subsetneq V_t$, i.e., missing some cut-yielding subproblems, may still lead to lower overall runtime. Adding fewer cuts at iteration~$t$ changes the master solution, which may cause a different set~$V_{t+1}$ at the next iteration---one that might yield stronger cuts than under $\Pi^*$. Hence, $T(\Pi^*)$ is a reference point rather than a proven lower bound on achievable subproblem solving time. We focus on approximating the oracle's behavior as a well-defined and computable benchmark.

\begin{remark}[Cut Prediction Oracle Bound]
    \label{rem:perfect-selection}
    Let $\tau_s^t$ denote the time to solve subproblem $s$ at iteration $t$. The \emph{cut prediction oracle bound} $T(\Pi^*)$ is the total subproblem solving time when solving only cut-yielding subproblems:
    \begin{equation}
        T(\Pi^*) = \sum_{t=1}^T \sum_{s \in V_t} \tau_s^t \leq \sum_{t=1}^T \sum_{s \in S} \tau_s^t = T(\Pi_{\text{all}}).
    \end{equation}
    Since the same set of cuts is added (i.e., the same scenarios are found infeasible with identical dual rays), the runtime behavior is equal otherwise, apart from effects of warmstarting. 
\end{remark}

In general, this cut prediction oracle is not computable a priori, but rather serves as a benchmark. Our goal is to design a scoring mechanism that approximates the oracle's performance by prioritizing subproblems likely to generate cuts.

\subsection{Subproblem Scoring}
\label{sec:subproblem-scoring}

For each scenario $s \in S$, we maintain a score $\texttt{score}(s)\in[0,1]$ that estimates its likelihood of generating a violated cut. We present a machine learning-based scoring approach that uses network topology features and historical performance metrics (Section~\ref{sec:ml-score}).

\subsubsection{Score Components}
\label{sec:score-components}

We utilize information derived from the master problem solution, in conjunction with the failure scenario and from historical scenario performance.
The choice of components is guided by two principles: (1)~historical cut generation metrics ($\rho_s$, $\tau_s$, $\bar{\nu}_s$, $\zeta_s$) capture temporal patterns observable in any Benders decomposition, analogous to multi-criteria cut scoring in MIP solvers~\cite{deza2023machine,turner2023cutting}; (2)~topology-aware features ($\kappa_e$, $\psi_e$, $\mu_e$, $\beta_e$) encode the structural importance of the failed edge in the current solution.
Let $(\bar{y}, \bar{f})$ denote the current master solution where $\bar{y}_{e,m}$ is the number of modules of type $m$ installed on link $e$, and $\bar{f}_{s_0,d,a}$ is the flow of demand $d$ on arc $a$ in the base case scenario $s_0$ (no failures).
Let $C_s$ denote the set of iterations where scenario $s$ generated a cut, and $I_s$ denote iterations where $s$ was examined.
Let $\nu_s^{t'}$ denote the violation magnitude of the Benders cut generated by scenario $s$ at iteration $t' \in C_s$, defined as
$\nu_s^{t'} = -\sum_{d \in D} \sum_{v \in V} \sigma_{d,v}^s b_{d,v} - \sum_{e \notin F_s} \pi_{e}^s \sum_{m \in M_e} u_{e,m} \bar{y}_{e,m},$
i.e., the amount by which the current solution violates the generated cut~\eqref{eq:benders_cut}.
Let betweenness centrality $\beta_e$ be the fraction of shortest paths (in the network topology, without considering modules) between each demand pair that contain edge $e$. Using these definitions, we introduce eight score components as described in Table~\ref{tab:score-components}.
Note that $\nu_s^{t'}$ for feasibility cuts is not normalized, since the dual ray may scale arbitrarily. Any observed effects might therefore either be due to an implicit normalization by the solver or due to the role of $\nu_s^{t'}$ as a proxy for the reliability.

\begin{table}[ht]
\centering

\caption{Score components for scenario prioritization.}
\label{tab:score-components}
\scalebox{0.70}{
\begin{tabular}{llp{4.9cm}}
\toprule
Name & Definition & Description \\
\midrule
Reliability & $\rho_s = |C_s| / |I_s|$ & Fraction of times $s$ produced a cut when solved \\
Total Share & $\tau_s = |C_s| / \sum_{s' \in S} |C_{s'}|$ & Fraction of all cuts produced by $s$ \\
Average Violation & $\bar{\nu}_s = \frac{1}{|C_s|} \sum_{t' \in C_s} \nu_s^{t'}$ & Mean violation magnitude when $s$ generated cuts \\
Staleness & $\zeta_s = t - \max\{t' : t' \in I_s\}$ & Iterations since $s$ was last examined \\
Failed Edge Capacity & $\kappa_e = \sum_{m \in M} u_{e,m} \bar{y}_{e,m}$ & Total capacity on failed edge $e \in F_s$ \\
Failed Edge Flow & $\psi_e = \sum_{d \in D} (|\bar{f}_{s_0,d,a}| + |\bar{f}_{s_0,d,a'}|)$ & Base case flow on failed edge $e \in F_s$ \\
Failed Edge Utilization & $\mu_e = \psi_e / \kappa_e$ & Flow-to-capacity ratio on $e\in F_s$ (when $\kappa_e > 0$) \\
Betweenness Centrality & $\beta_e$ & Topological centrality of edge $e \in F_s$ (computed once) \\
\bottomrule
\end{tabular}
}
\end{table}

A straightforward approach would combine score components via weighted linear combination $R(s) = w_\rho \rho_s + w_\tau \tau_s + w_\zeta \zeta_s + \ldots$ with min-max normalization. However, determining appropriate weights is challenging: they vary across instances and solution phases. For problems with long solve times, static weights become problematic as the optimal balance shifts dynamically. This motivates the machine learning approach described next, which learns feature importance automatically.

\subsubsection{Machine Learning-Based Scoring}
\label{sec:ml-score}

Our regression-based approach (in the following, referred to as machine learning, or ``ML'') predicts scenario infeasibility using logistic regression trained online during the Benders decomposition process. We construct an 8-dimensional feature vector $\phi(s) \in \R^{8}$ for each scenario $s$ combining the failed edge capacity $\kappa_e$, flow $\psi_e$, utilization $\mu_e$, betweenness centrality $\beta_e$, and the historical components $\bar{\nu}_s$, $\rho_s$, $\tau_s$, $\zeta_s$ described above. Features are normalized to zero mean and unit variance using Welford's online algorithm~\cite{welford1962note}, with normalization parameters accumulated across the entire solution process (not reset during stabilization rounds).

\paragraph*{Logistic Regression Model}

We use logistic regression for its low computational complexity during online training and superior interpretability compared to more sophisticated models like neural networks or ensemble methods.

Logistic regression models the probability of binary outcomes via a sigmoid transformation of a linear predictor. The model predicts infeasibility probability as:
{
\begin{equation}
    p(s \text{ infeasible} \mid \phi(s); w, b) = \sigma(w^\top \phi(s) + b)
\end{equation}
}
where $w \in \R^{8}$ is the weight vector, $b \in \R$ is the bias term, and $\sigma(z) = 1/(1+e^{-z})$ is the logistic function mapping $\R \to (0,1)$.

This formulation allows the model to learn non-uniform importance of features: features strongly predictive of infeasibility receive large absolute weights, while less informative features receive weights near zero. The sigmoid ensures output values are valid probabilities.

\paragraph*{Online Training via Gradient Descent}

After solving the subproblem for scenario $s$ at iteration $t$, we observe the outcome $z_s^t \in \{0,1\}$ where $z_s^t = 1$ if the subproblem was infeasible (generated a cut) and $z_s^t = 0$ if feasible. We update the model parameters to minimize the logistic loss:
{
\begin{equation}
    \ell(w, b) = -\log p(z_s^t \mid \phi(s); w, b)
\end{equation}
}
\noindent
Applying gradient descent with L2 regularization yields the update rules:
{
\begin{align}
    w &\gets w - \alpha \left[ (\hat{z}_s^t - z_s^t) \phi(s) + \lambda w \right] \label{eq:sgd_w}\\
    b &\gets b - \alpha (\hat{z}_s^t - z_s^t) \label{eq:sgd_b}
\end{align}
}
where $\hat{z}_s^t = \sigma(w^\top \phi(s) + b)$ is the predicted probability, $\alpha > 0$ is the learning rate controlling step size, and $\lambda \geq 0$ is the regularization parameter preventing overfitting by penalizing large weights.

The gradient term $(\hat{z}_s^t - z_s^t)$ represents the prediction error: positive when the model overestimates infeasibility probability, negative when underestimating. The regularization term $\lambda w$ in~\eqref{eq:sgd_w} shrinks weights toward zero, favoring simpler models that generalize better.
Note that these online updates yield an approximation of the model that would result from batch retraining on all accumulated data. We choose online updates for their $O(1)$ per-observation cost, which keeps the training overhead negligible; investigating the potential benefit of periodic batch retraining on small instances is left for future work.
\subsection{Subproblem Selection}
\label{sec:partial-subproblem-solving}

\begin{table}[ht]
\centering
\caption{Stopping criteria for partial subproblem solving and naming convention. One or more can be imposed, and a Benders iteration stops subproblem solving when any criterion is met. Stabilization rounds solve all subproblems to prevent selection bias and provide additional training data.}
\label{tab:stopping-criteria}
\scalebox{0.72}{
\begin{tabular}{llp{10cm}}
\toprule
Name & Identifier & Description \\
\midrule
Cut limit & \texttt{$k$C} & Stop if $k$ cuts have been generated \\
Relative solve limit & \texttt{$p$P} & Stop after examining a proportion $p \in [0,1]$ of all scenarios, i.e., $\lceil p \cdot |S| \rceil$ subproblems \\
Consecutive misses & \texttt{$m$M} & Stop if no cuts found for $m$ consecutive scenarios \\
Score limit & \texttt{$r$T} & Stop when the score of the next scenario falls below threshold $r$ \\
Time limit & --- & Stop if iteration time limit reached \\
\midrule
\multicolumn{3}{l}{\textbf{Stabilization} (solving all subproblems)} \\
\midrule
Stabilization & \texttt{$n$S} & Every $n$ iterations \\
Root node & \texttt{$n$R} & First $n$ rounds at root \\
\bottomrule
\end{tabular}}
\end{table}

We solve subproblems in order of decreasing score.
Using the stopping criteria introduced in \autoref{tab:stopping-criteria}, we aim to examine all cut-generating scenarios, and to stop once we are unlikely to find more cuts.
Since limited examination may miss some cut-generating scenarios and, over time, introduce a selection bias, we periodically perform \emph{stabilization rounds} where all scenarios are examined (disabling stopping criteria). This prevents systematic exclusion of subproblems and provides additional training data. 
Since initial training data is especially important, we always initialize the scores with at least one score initialization round, but up to $n$ stabilization rounds at the root node (\texttt{$n$R}).
At stabilization rounds, we reset the score components $\rho_s$, $\tau_s$, $\bar\nu_s$ and $\zeta_s$ to avoid bias from outdated historical data.

\section{Computational Experiments}
\label{sec:experiments}

\subsection{Experimental Setup}

We implement our Branch-and-Benders-cut framework with adaptive cut selection using Gurobi 12.0.3 \cite{gurobi}, Julia 1.12.4 \cite{bezanson2017julia} with JuMP 1.12 \cite{dunning2017jump}, and lazy constraint callbacks. To improve efficiency, we maintain a single subproblem template adjusted for each scenario rather than constructing new subproblems from scratch. We reuse the Gurobi environment, allowing for LP warmstarting. We do not perform separation on fractional solutions so as to isolate the impact of subproblem selection.

\autoref{tab:config} presents the configuration parameters for all evaluated methods.
We parameterize those configurations following an analysis of stopping criteria: For each (threshold, proportion, cut limit) individually, we examine how many cuts would have been generated under this setting. A distribution boxplot can be found in \autoref{fig:ml-train-parameter-boxplots} in \autoref{sec:appendix-subproblem-stats}. Based on this analysis, we select more conservative (calling more subproblems, possibly solving more unrequired ones) parameters as well as more aggressive ones to evaluate the impact of stopping criteria. In particular, setting \texttt{ML-0.1T-0.6P-10C-20S-5R} uses all three main stopping criteria, and sets the limit at the upper quartile of required cuts per iteration (see \autoref{fig:ml-train-parameter-boxplots} in \autoref{sec:appendix-subproblem-stats}). That way, we expect to capture most cut-generating scenarios while avoiding excessive unrequired subproblem solves.

\begin{table}[ht]
\centering
\caption{Configuration parameters for tested methods. ML methods use online logistic regression. See \autoref{tab:stopping-criteria} for parameter descriptions.}
\label{tab:config}
\resizebox{\textwidth}{!}{
\begin{tabular}{llccccc}
\toprule
Method & Selection & Threshold & Proportion & Cut Limit & Stabilization & Root Node Stabilization \\
\midrule
Standard & All & --- & --- & --- & --- & --- \\
\texttt{Random-0.5P-20S-5R} & Random & --- & 0.5 & --- & 20 & 5 rounds \\
\texttt{ML-0.1T-0.6P-10C-20S-5R} & Threshold+Prop & 0.1 & 0.6 & 10C & 20 & 5 rounds \\
\texttt{ML-0.5P-20S-5R} & Proportion & --- & 0.5 & --- & 20 & 5 rounds \\
\texttt{ML-0.5T-5C-20S-5R} & Threshold & 0.5 & --- & 5C & 20 & 5 rounds \\
\texttt{ML-0.7P-20S-5R} & Proportion & --- & 0.7 & --- & 20 & 5 rounds \\
\bottomrule
\end{tabular}
}
\end{table}

Out of the original SNDlib instances \cite{IdzikowskiOrlowskiRaackWoesnerWolisz2010b,KosterKutschkaRaack2010}, only 11 solve within our time limit of three hours, which yields an insufficient sample for meaningful evaluation.
We generate test instances from SNDlib benchmark networks by combining 2--5 subgraphs: (1) extract connected subgraphs via BFS from high-degree nodes using proportion $p \in [0.45, 0.65]$ to obtain $p \cdot |V|$ nodes per network; (2) merge subgraphs by connecting high-degree nodes; (3) generate N-1 failure scenarios ($|S| = |E|$); (4) filter scenarios that disconnect demand nodes. We generate 135 test instances with 16--229 scenarios, 8--57 nodes, 13--82 edges, and 18--857 demands. The online ML training overhead is negligible, comparable to solving a single subproblem per iteration. For stabilization frequency, we tested values of 5, 10, 20, 50, 100, and 200 iterations, with 20 providing the best balance for SND problems.

All experiments use machines with 2x Intel Xeon L5630 Quad Core processors at 2.13 GHz (8 cores, 16GB RAM). A three-hour time limit is imposed per instance. Statistical comparisons use the shifted geometric mean (shift of 10) to handle zero values.

The parameters for the regression model were determined using a grid search on a set of 14 instances, generated synthetically by the above procedure. We tried configurations with learning rates $\alpha \in \{0.02, 0.05, 0.075, 0.1, 0.15\}$ and regularization $\lambda \in \{0.05, 0.02, 0.01, 0.005, 0.001\}$. The selected parameters ($\alpha=0.075$, $\lambda=0.02$) provided the best runtime performance on the tuning set.
\subsection{Results}

\begin{figure}[htbp]
    \centering
    \includegraphics[width=\textwidth]{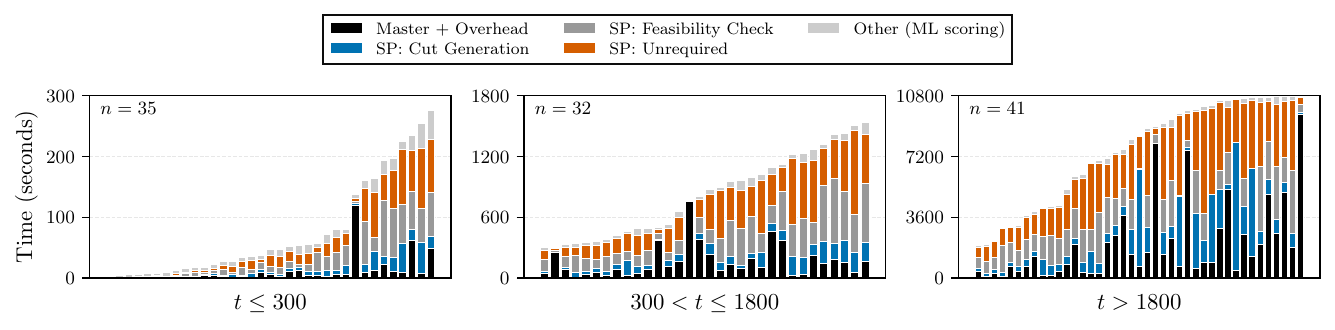}
    \caption{Time breakdown across $n=108$ instances. Master time (black), cut-generating subproblems (blue), and feasibility-proving subproblems in cut-free iterations (grey) represent necessary computation; unrequired subproblems (orange) represent wasted computation.}
    \label{fig:time-comparison-standard-oracle}
\end{figure}

The potential for subproblem selection is established through the cut prediction oracle bound (see \autoref{rem:perfect-selection}). As illustrated in \autoref{fig:time-comparison-standard-oracle} (and described in detail in \autoref{tab:subproblem-stats} in \autoref{sec:appendix-subproblem-stats}), 34.4\% of total solve time in standard Benders is spent on unrequired solving of subproblems. Since solving all subproblems that yield cuts does not always lead to the best runtime, this is not necessarily an upper bound on the achievable speedup. However, we consider it as a dimension for comparison.

\paragraph*{Runtime Performance}

\begin{figure}[htbp]
    \centering
    \includegraphics[width=0.75\textwidth]{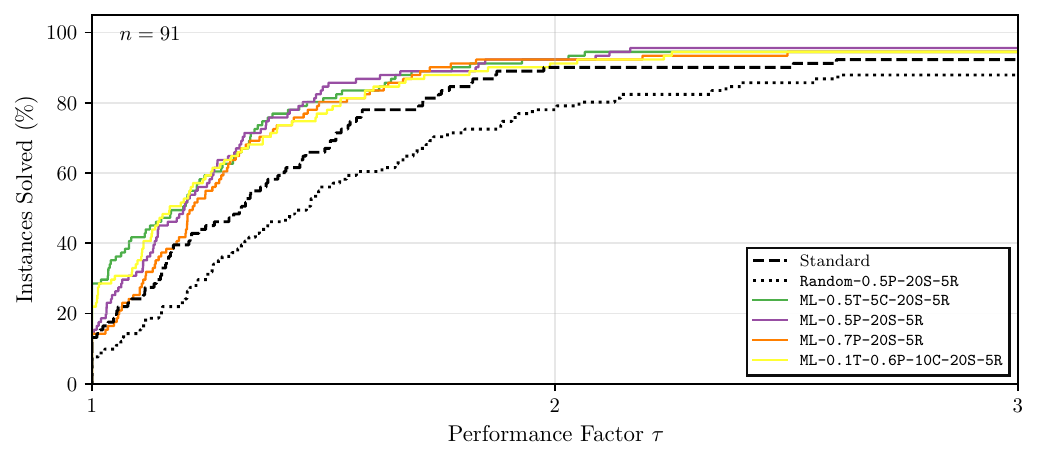}
    \caption{Performance profile on all $n=91$ instances solved by at least one method. The performance factor $\tau$ on the $x$-axis is the ratio of a method's runtime to the best method's runtime on each instance; $\rho(\tau)$ on the $y$-axis is the fraction of instances where the method is within factor $\tau$ of the best. Random selection performs worse than standard Benders, regression-based methods achieve modest gains.}
    \label{fig:performance-profile}
\end{figure}

\begin{table}[t]
\centering
\caption{Summary of runtime performance across all methods ($n=79$ instances solved by all methods). Solved: instances solved to optimality; Sum/Mean/Sh. Geom. Mean$^\dagger$: computed over instances solved by all methods. The shifted geometric mean (with shift $s=10$) is $\bigl(\prod_{i=1}^{n}(x_i + s)\bigr)^{1/n} - s$, which dampens the influence of very small values. Wilcoxon signed-rank test~\cite{wilcoxon1945individual} (one-sided): a nonparametric test for matched pairs that tests if each method is significantly faster than Standard; $p < 0.05$ indicates the observed improvement is unlikely under the null hypothesis of no difference. Significance levels: * $p<0.05$, ** $p<0.01$, *** $p<0.001$, n.s. = not significant.}
\label{tab:methods-summary}
\scalebox{0.75}{
\begin{tabular}{lrrrrrr}
\toprule
Method & Solved & Sum$^\dagger$ & Mean$^\dagger$ & Sh. Geom. & Wilcoxon & Signif- \\
 &  &  &  & Mean$^\dagger$ & $p$ & icance \\
\midrule
Standard & 86 & 100941.2 & 1277.7 & 300.5 & --- & --- \\
\texttt{Random-0.5P-20S-5R} & 82 & 112118.3 & 1419.2 & 334.8 & 0.9996 & n.s. \\
\texttt{ML-0.5T-5C-20S-5R} & 86 & 86715.2 & 1097.7 & 270.8 & 0.0130 & * \\
\texttt{ML-0.5P-20S-5R} & 89 & 96097.3 & 1216.4 & 285.2 & 0.0136 & * \\
\texttt{ML-0.7P-20S-5R} & 88 & 100536.4 & 1272.6 & 288.3 & 0.2080 & n.s. \\
\texttt{0.1T-0.6P-10C-20S-5R} & 87 & 92977.6 & 1176.9 & 275.3 & 0.0081 & ** \\
\bottomrule
\end{tabular}
}
\end{table}

As illustrated in \autoref{fig:performance-profile} and detailed in \autoref{tab:methods-summary}, three out of four tested ML configurations achieve statistically significant speedups over standard Benders, the best being \texttt{ML-0.1T-0.6P-10C-20S-5R}. The reduction from 300.5 to 270.8 represents a speedup of 9.9\%, as opposed to a ``theoretical'' maximum of 34.4\% (see \autoref{tab:subproblem-stats} in \autoref{sec:appendix-subproblem-stats}). Since we aim for reducing unrequired subproblem solving times, the performance profile reduced to only those times shows a more pronounced improvement (see \autoref{fig:adjusted-sp-performance-profile} in \autoref{sec:appendix-sp-ratio-speedup}).
The full results can be found in \autoref{tab:results} in \autoref{sec:appendix-results}.
We observe that our scoring method indeed outperforms random selection, which is statistically significantly worse than standard Benders (one-sided Wilcoxon signed-rank test, $p=0.9996$, indicating random is slower). This confirms that naive subproblem selection strategies can degrade performance, and that our proposed scoring mechanism is able to prioritize important subproblems.

\paragraph*{Convergence Quality}

Beyond runtime performance, primal-dual integral analysis (detailed results in \autoref{tab:pd-integrals} in \autoref{sec:appendix-pd-integrals}) evaluates convergence quality. All four ML configurations achieve statistically significant improvements in convergence quality over standard Benders, with the best method \texttt{ML-0.5T-5C-20S-5R} reducing the shifted geometric mean from 1336.3 to 1145.9 billion (14.3\% improvement). Random selection shows no significant difference from standard Benders ($p=0.12$). These results demonstrate that ML-based subproblem selection not only reduces solve time but also improves the quality of the search by better exploring the solution space.

The ML model achieves an average accuracy (fraction of correct predictions) of 84.4\%, precision (true positives among predicted positives) of 33.4\%, recall (true positives among actual positives) of 70.9\%, and F1-score (harmonic mean of precision and recall) of 45.2\% (see \autoref{tab:ml_aggregated} in \autoref{sec:appendix-ml-analysis}). The mediocre scores mean that we are ``missing'' cuts, which can lead to worse algorithm convergence that can also be seen in the runtime table (\autoref{tab:results} in \autoref{sec:appendix-results}).

The ML model's imperfect predictions can be attributed to two main factors: first, the features used may not fully capture the complex relationships around scenario violations (as can be observed by the high variance in feature importance, see \autoref{fig:feature-importance} in \autoref{sec:appendix-ml-analysis}) and second, the online training approach may not obtain sufficient training data for robust learning, especially when instances solve quickly.

To counteract the former, we also tested locality-based features (e.g., utilization in the $n$-hop neighborhood of failed edges) but they did not receive significant weights in our preliminary experiments. Addressing the latter, we tried pre-training on similar instances, but as \autoref{fig:cut-yield-evolution} in \autoref{sec:appendix-subproblem-stats} shows, cut generation patterns evolve substantially during solving, and so do the weights of the ML model. Additionally, we implemented root node stabilization that serves the purpose of establishing a strong initial ML model and facilitating early primal heuristics. However, further research is needed to develop more reliable prediction models that can consistently deliver robust performance improvements across diverse instances.

\section{Conclusions and Future Work}
\label{sec:conclusion}

This paper investigated adaptive subproblem selection in Benders decomposition, motivated by the empirical observation that 52.1\% of subproblems in our variant of the survivable network design problem unnecessarily solved: The rest either yield cuts (5.6\%), or verify feasibility in cut-free iterations (42.4\%). We developed a multi-criteria scoring mechanism combining historical performance metrics with topology-aware features, trained online via logistic regression. Multiple stopping criteria (score thresholds, proportional limits, cut limits) adaptively determine how many subproblems to solve per iteration. The online training overhead is negligible---comparable to solving a single subproblem---making the approach practical for real-time deployment.

Computational experiments on 135 survivable network design instances reveal that in standard Benders with full subproblem enumeration, 34.4\% of total solve time is spent on unrequired subproblem solving.
Our best runtime configuration (\texttt{ML\mbox{-}0.1T\mbox{-}0.6P\mbox{-}10C\mbox{-}20S\mbox{-}5R}) achieves 9.9\% speedup over standard Benders (270.8s vs 300.5s shifted geometric mean, $p=0.0081$), while the best convergence configuration (\texttt{ML\mbox{-}0.5T\mbox{-}5C\mbox{-}20S\mbox{-}5R}) reduces primal-dual integrals by 14.3\% ($p<0.001$). Three of four ML configurations achieve statistically significant runtime improvements, and all achieve statistically significant improvements in primal-dual integrals. While the scoring performance is modest (the regression model achieves only 33.4\% precision, 70.9\% recall, and 45.2\% F1), these results establish that adaptive subproblem selection based on learned importance can improve Benders decomposition performance. Combining multiple stopping criteria appears to be a promising direction, as our best method employs three types simultaneously. Random selection degrades performance significantly ($p=0.9996$), confirming that informed prioritization is necessary.

\paragraph*{Future Research Directions}
Multiple future research directions can be identified: First, more sophisticated prediction models---such as graph neural networks---could better capture network topology and temporal evolution of cut generation patterns, potentially improving precision and allowing for preconditioning on a set of similar instances. Second, rather than using fixed thresholds, dynamic adaptation of stopping criteria during solving could adjust exploration-exploitation tradeoffs based on observed performance, also taking into account cut usefulness. Third, analyzing feature importance across different solving stages could reveal which metrics are most predictive early versus late in the algorithm, enabling stage-specific scoring; similarly, restricting the model to a subset of high-importance features may improve robustness and interpretability. Fourth, a hybrid approach could apply standard Benders decomposition (possibly with separation on fractional solutions) in early iterations and switch to ML-based selection in later iterations where cut yield declines. Finally, auxiliary (problem-specific) optimization techniques based on min-cut computations could provide tighter bounds on scenario criticality, complementing the learned scoring mechanism. These directions could help close the gap between current performance and the theoretical potential identified through our analysis of subproblem selection potential.

\subsubsection*{Acknowledgements}
I thank Oliver Gaul for his help in previous work on adaptive subproblem selection and for insightful discussions, Marco L\"{u}bbecke for guidance and review, and Alexander Helber for helpful comments.
I also thank the anonymous reviewers for their constructive feedback, which improved the presentation of this paper.

\subsubsection*{Code Availability}
Source code is available at \url{https://github.com/tidonk/BendersNetworkDesign.jl}.


\bibliography{references}

\newpage
\appendix
\section*{Appendix}

\section{Subproblem Statistics}
\label{sec:appendix-subproblem-stats}

This section provides detailed subproblem statistics for all instances and the distribution of stopping criteria parameter values.

\begin{table}[ht]
\centering
\caption{Time and count breakdown for ML-train method across instances that completed within the time limit ($n=86$). Unrequired subproblems are those that do not yield cuts in iterations where other subproblems do. Cut Generation shows subproblems that yielded cuts. Feasibility Proving shows subproblems in cut-free iterations (required to verify feasibility). Other includes ML scoring and callback overhead.}
\label{tab:subproblem-stats}
\scalebox{0.85}{
\begin{tabular}{lrrrrr}
\toprule
& \multicolumn{3}{c}{Time} & \multicolumn{2}{c}{Subproblems} \\
\cmidrule(lr){2-4} \cmidrule(lr){5-6}
Category & Time (s) & \% Total & \% SP & Count & \% SP \\
\midrule
Master Problem & 87,626 & 25.8\% & --- & --- & --- \\
SP: Cut Generation & 56,712 & 16.7\% & 23.5\% & 20,485 & 5.6\% \\
SP: Feasibility Proving & 68,367 & 20.2\% & 28.3\% & 155,251 & 42.4\% \\
SP: Unrequired & \textbf{116,594} & \textbf{34.4\%} & \textbf{48.2\%} & \textbf{190,788} & \textbf{52.1\%} \\
Other (ML, overhead) & 9,786 & 2.9\% & --- & --- & --- \\
\midrule
Total & 339,086 & 100.0\% & --- & 366,524 & 100.0\% \\
\bottomrule
\end{tabular}
}
\end{table}

\begin{figure}[htbp!]
    \centering
    \includegraphics[width=\textwidth]{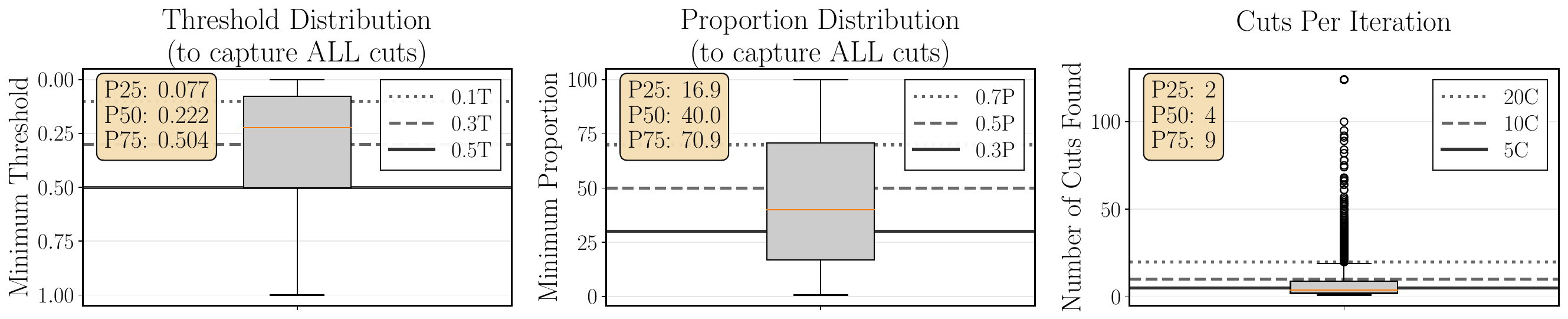}
    \caption{Distribution of minimum parameter values required to capture all cuts per iteration during ML training, with reference lines showing tested configuration values. Box plots show the median (orange line), interquartile range (box), whiskers extending to $1.5\times$ the interquartile range, and the mean (white diamond); points beyond whiskers are outliers.}
    \label{fig:ml-train-parameter-boxplots}
\end{figure}

\begin{figure}[htbp!]
    \centering
    \includegraphics[width=0.65\textwidth]{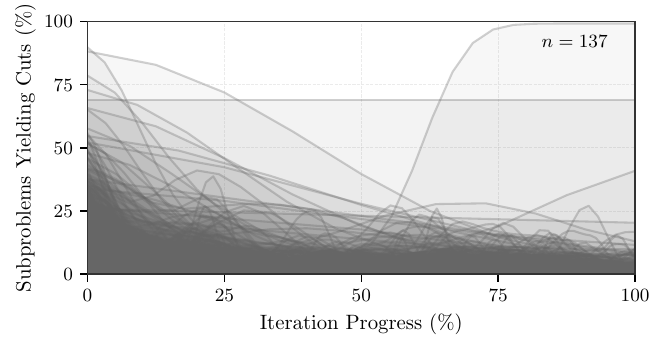}
    \caption{Cut yield evolution across iterations for full enumeration, showing the fraction of examined subproblems that yield cuts. Patterns decline from 10--40\% early on to 0--10\% near convergence.}
    \label{fig:cut-yield-evolution}
\end{figure}

\clearpage
\section{Regression Model Analysis}
\label{sec:appendix-ml-analysis}

Classification metrics for the logistic regression model, aggregated across all instances. Accuracy is the fraction of correct predictions. Precision is the fraction of predicted cut-yielding subproblems that actually yield cuts (true positives / predicted positives). Recall is the fraction of actual cut-yielding subproblems that the model correctly identifies (true positives / actual positives). F1-score is the harmonic mean of precision and recall, $\text{F1} = 2 \cdot \text{precision} \cdot \text{recall} / (\text{precision} + \text{recall})$.

\begin{table}[ht]
\centering
\caption{Aggregated ML Model Classification Metrics (\%) by Configuration. Values show geometric mean across all instances. Standard includes full enumeration, enabling slightly better performance due to more training data.}
\label{tab:ml_aggregated}
\scalebox{0.75}{
\begin{tabular}{lrrrr}
\toprule
Configuration & Accuracy & Precision & Recall & F1-Score \\
\midrule
Standard & 86.0 & 34.3 & 71.0 & 46.1 \\
\texttt{ML-0.5P-20S-5R} & 83.8 & 32.8 & 70.4 & 44.6 \\
\texttt{ML-0.5T-5C-20S-5R} & 83.9 & 32.6 & 69.8 & 44.3 \\
\texttt{ML-0.7P-20S-5R} & 84.9 & 34.0 & 71.5 & 46.0 \\
\texttt{ML-0.1T-0.6P-10C-20S-5R} & 83.5 & 33.2 & 71.7 & 45.2 \\
\midrule
Average & 84.4 & 33.4 & 70.9 & 45.2 \\
\bottomrule
\end{tabular}
}
\end{table}

\begin{figure}[htbp!]
    \centering
    \includegraphics[width=0.7\textwidth]{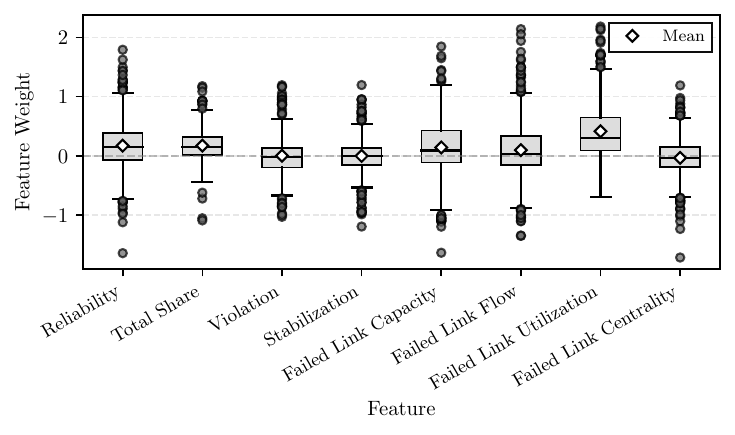}
    \caption{Distribution of learned feature weights across instances. Feature reference in \autoref{tab:score-components}. Box plots show the median (line), interquartile range (box), whiskers extending to $1.5\times$ the interquartile range, and the mean (white diamond); points beyond whiskers are outliers.}
    \label{fig:feature-importance}
\end{figure}

\clearpage

\section{Subproblem Selection Quality}
\label{sec:appendix-sp-ratio-speedup}

Hit quality measures the efficiency of subproblem selection as the ratio of cuts found to subproblems examined per iteration, averaged across all iterations. The adjusted subproblem performance profile isolates the effect of selection by deducting the necessary subproblem time (oracle baseline) from each method's total subproblem time.

\begin{figure}[htbp!]
    \centering
    \includegraphics[width=0.7\textwidth]{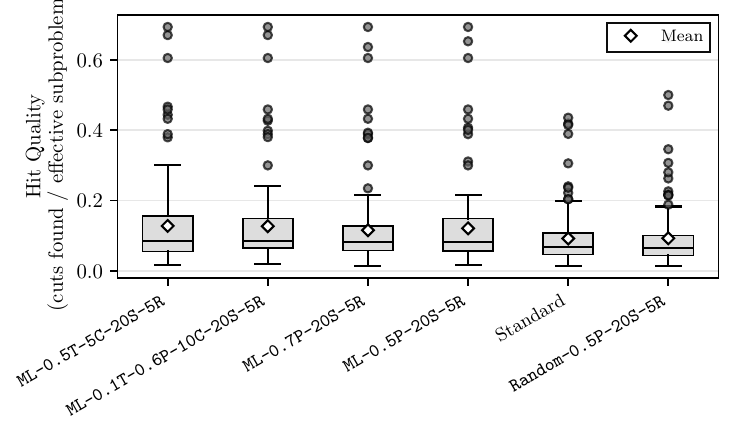}
    \caption{Distribution of hit quality (cuts found / subproblems examined) across methods. Higher values indicate more efficient subproblem selection. Box plots show the median (line), interquartile range (box), whiskers extending to $1.5\times$ the interquartile range, and the mean (white diamond); points beyond whiskers are outliers.}
    \label{fig:hitquality-boxplot}
\end{figure}

\begin{figure}[htbp!]
    \centering
    \includegraphics[width=0.7\textwidth]{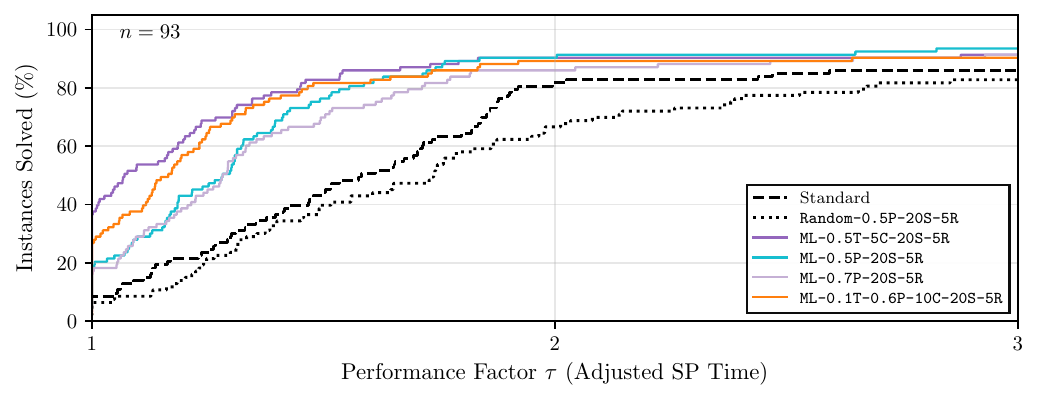}
    \caption{Performance profile based on adjusted subproblem time, which deducts the necessary subproblem time (oracle baseline) from each method's total subproblem time.}
    \label{fig:adjusted-sp-performance-profile}
\end{figure}

\clearpage

\section{Detailed Solve Times}
\label{sec:appendix-results}

\begin{table}[htbp!]
\centering
\caption{Solve times (seconds) for all instances. Timeouts: optimality gaps in brackets. Bold: methods faster than Standard (if Standard solved), or methods that solved (if Standard did not). Daggers ($^\dagger$): all methods optimal. Sum/Mean/Sh.\ Geom.\ Mean$^\dagger$: over instances solved by all. Wilcoxon test (one-sided): solved use time; unsolved use time-limit+gap\%. Significance: * $p<0.05$, ** $p<0.01$, *** $p<0.001$. Summary statistics (bottom right) are computed over both columns.}
\label{tab:results}
\scalebox{0.5}{
\begin{tabular}{@{}lr@{\,}r@{\,}r@{\,}r@{\,}r@{\,}r@{\quad}|@{\quad}lr@{\,}r@{\,}r@{\,}r@{\,}r@{\,}r@{}}
\toprule
 & \multicolumn{2}{c}{\textbf{Baseline}} & \multicolumn{4}{c}{\textbf{ML-Based}} & & \multicolumn{2}{c}{\textbf{Baseline}} & \multicolumn{4}{c}{\textbf{ML-Based}} \\
\cmidrule(lr){2-3} \cmidrule(lr{2em}){4-7} \cmidrule(lr){9-10} \cmidrule(lr){11-14}
Instance & \makecell{Standard} & \makecell{\texttt{Rand.}\\\texttt{0.5P}} & \makecell{\texttt{ML}\\\texttt{0.5T-5C}} & \makecell{\texttt{ML}\\\texttt{0.5P}} & \makecell{\texttt{ML}\\\texttt{0.7P}} & \makecell{\texttt{ML}\\\texttt{0.1T-0.6P}\\\texttt{10C}} & Instance & \makecell{Standard} & \makecell{\texttt{Rand.}\\\texttt{0.5P}} & \makecell{\texttt{ML}\\\texttt{0.5T-5C}} & \makecell{\texttt{ML}\\\texttt{0.5P}} & \makecell{\texttt{ML}\\\texttt{0.7P}} & \makecell{\texttt{ML}\\\texttt{0.1T-0.6P}\\\texttt{10C}} \\
\midrule
\texttt{abilene}$^\dagger$ & 6.9 & 8.5 & \textbf{5.1} & \textbf{5.2} & \textbf{5.0} & \textbf{4.8} & \multicolumn{7}{c@{}}{\scalebox{0.6}{$\vdots$}} \\
\texttt{atlanta}$^\dagger$ & 56.2 & 89.2 & \textbf{55.6} & 68.3 & \textbf{53.7} & \textbf{54.3} & \texttt{075\_n5s142} & (0.2\%) & (5.4\%) & \textbf{(0.1\%)} & (0.6\%) & (1.8\%) & (0.8\%) \\
\texttt{di-yuan}$^\dagger$ & \textbf{735.7} & 1743.3 & 765.6 & 961.8 & 814.9 & 1091.1 & \texttt{076\_n5s122} & (3.4\%) & (9.5\%) & \textbf{(0.6\%)} & (12.5\%) & \textbf{(2.2\%)} & (10.2\%) \\
\texttt{france}$^\dagger$ & 248.7 & \textbf{205.1} & \textbf{175.9} & \textbf{210.8} & \textbf{162.8} & \textbf{183.6} & \texttt{081\_n2s160} & (76.3\%) & \textbf{(72.2\%)} & \textbf{(73.0\%)} & \textbf{(69.6\%)} & \textbf{(67.3\%)} & \textbf{(72.9\%)} \\
\texttt{geant}$^\dagger$ & 614.3 & 1581.4 & \textbf{337.5} & \textbf{428.9} & \textbf{359.2} & \textbf{370.7} & \texttt{083\_n2s76} & (1.7\%) & (7.8\%) & \textbf{(0.1\%)} & (16.1\%) & (5.8\%) & (8.3\%) \\
\texttt{newyork}$^\dagger$ & 841.4 & 1134.5 & 1103.0 & 983.9 & \textbf{822.0} & 1080.7 & \texttt{084\_n2s38} & \textbf{4829.9} & (0.1\%) & (0.2\%) & 5518.6 & (0.3\%) & 8962.0 \\
\texttt{nobel-eu}$^\dagger$ & 2419.8 & \textbf{2157.9} & \textbf{1418.7} & \textbf{1616.4} & \textbf{1797.9} & \textbf{1478.2} & \texttt{085\_n2s83}$^\dagger$ & 600.6 & 671.5 & 728.5 & \textbf{563.1} & 674.7 & \textbf{557.3} \\
\texttt{nobel-germany}$^\dagger$ & 40.9 & \textbf{36.8} & \textbf{36.4} & \textbf{38.4} & 43.0 & \textbf{28.9} & \texttt{086\_n2s66}$^\dagger$ & 1838.4 & 2670.9 & \textbf{1710.8} & \textbf{1480.3} & \textbf{1511.3} & \textbf{1418.4} \\
\texttt{pdh}$^\dagger$ & \textbf{103.8} & 178.8 & 214.4 & 154.1 & 126.8 & 143.9 & \texttt{087\_n2s86} & (41.9\%) & (45.5\%) & \textbf{(29.0\%)} & \textbf{(38.9\%)} & \textbf{(39.3\%)} & \textbf{(28.7\%)} \\
\texttt{polska}$^\dagger$ & \textbf{6.5} & 6.9 & 7.8 & 7.3 & 7.1 & 6.8 & \texttt{088\_n2s55}$^\dagger$ & 302.5 & 377.1 & 331.0 & \textbf{200.7} & 319.9 & \textbf{199.9} \\
\texttt{sun}$^\dagger$ & 4446.0 & \textbf{2031.0} & \textbf{1612.6} & \textbf{1379.3} & \textbf{2022.1} & \textbf{1644.3} & \texttt{091\_n3s62}$^\dagger$ & 890.6 & 960.3 & 1060.6 & 909.7 & 955.9 & \textbf{860.4} \\
\texttt{001\_n2s26}$^\dagger$ & 7.4 & 7.7 & 8.3 & \textbf{7.3} & \textbf{6.7} & \textbf{6.1} & \texttt{092\_n3s26}$^\dagger$ & \textbf{3.0} & 3.6 & 3.7 & 3.7 & 4.4 & 4.2 \\
\texttt{003\_n2s76}$^\dagger$ & 2754.0 & 3537.4 & 4296.5 & \textbf{2417.3} & 2779.1 & 3716.2 & \texttt{093\_n3s90}$^\dagger$ & 4849.0 & \textbf{4714.0} & \textbf{2453.8} & \textbf{3245.9} & 5371.7 & 5485.4 \\
\texttt{004\_n2s33}$^\dagger$ & 36.7 & \textbf{26.0} & \textbf{25.3} & \textbf{31.1} & \textbf{26.5} & 40.2 & \texttt{094\_n3s69}$^\dagger$ & 500.7 & \textbf{438.0} & \textbf{427.8} & \textbf{462.8} & \textbf{285.3} & \textbf{364.4} \\
\texttt{006\_n2s51}$^\dagger$ & 510.0 & 560.5 & 554.8 & \textbf{414.2} & \textbf{468.7} & 513.5 & \texttt{095\_n3s55}$^\dagger$ & 276.4 & \textbf{240.2} & 298.6 & 1408.2 & 288.8 & 288.0 \\
\texttt{007\_n2s43}$^\dagger$ & 1196.8 & 1455.6 & \textbf{1075.8} & \textbf{1075.9} & 1269.4 & \textbf{1055.0} & \texttt{098\_n3s94} & 7637.6 & (2.9\%) & 9625.8 & \textbf{6300.3} & 8034.4 & 9498.9 \\
\texttt{008\_n2s156} & \textbf{(51.9\%)} & (54.0\%) & (56.9\%) & (57.1\%) & (53.4\%) & (56.3\%) & \texttt{099\_n3s102} & \textbf{(31.6\%)} & (33.0\%) & (36.4\%) & (36.4\%) & (36.4\%) & (36.4\%) \\
\texttt{009\_n2s35}$^\dagger$ & 12.6 & 17.7 & 13.1 & \textbf{11.2} & \textbf{11.7} & \textbf{10.9} & \texttt{104\_n4s152} & (28.5\%) & (31.6\%) & (34.7\%) & (34.4\%) & (30.3\%) & \textbf{(28.5\%)} \\
\texttt{010\_n2s42}$^\dagger$ & 293.1 & 295.2 & \textbf{190.5} & \textbf{282.0} & 768.1 & \textbf{254.8} & \texttt{112\_n5s78}$^\dagger$ & 653.4 & 730.7 & \textbf{481.4} & \textbf{627.8} & \textbf{600.1} & \textbf{431.3} \\
\texttt{011\_n3s64}$^\dagger$ & 681.8 & 738.0 & \textbf{580.0} & \textbf{495.3} & \textbf{582.8} & \textbf{637.0} & \texttt{113\_n5s131} & (1.7\%) & (5.6\%) & \textbf{(1.0\%)} & \textbf{(1.2\%)} & (1.9\%) & \textbf{(1.4\%)} \\
\texttt{012\_n3s78}$^\dagger$ & 3871.1 & 4471.7 & 4282.4 & \textbf{3339.7} & 5584.7 & \textbf{3770.7} & \texttt{115\_n5s150} & 8807.6 & 10332.3 & \textbf{8405.6} & \textbf{8528.3} & 9563.1 & (0.1\%) \\
\texttt{013\_n3s87} & 4774.8 & (0.6\%) & (0.2\%) & \textbf{3621.0} & \textbf{4409.6} & 5823.2 & \texttt{121\_n2s27}$^\dagger$ & 51.1 & \textbf{43.6} & 63.3 & \textbf{32.8} & 58.2 & \textbf{49.9} \\
\texttt{014\_n3s108} & (53.1\%) & (53.6\%) & \textbf{(47.7\%)} & \textbf{(52.5\%)} & \textbf{(48.6\%)} & (56.7\%) & \texttt{122\_n2s54}$^\dagger$ & 181.2 & \textbf{164.0} & \textbf{143.3} & \textbf{152.7} & \textbf{151.7} & \textbf{145.2} \\
\texttt{015\_n3s91}$^\dagger$ & 959.3 & \textbf{703.3} & \textbf{944.7} & \textbf{887.9} & \textbf{848.3} & 963.0 & \texttt{123\_n2s90}$^\dagger$ & 1021.8 & 1282.3 & 1049.3 & 1045.7 & 1092.8 & \textbf{941.6} \\
\texttt{017\_n3s84}$^\dagger$ & 7273.0 & \textbf{5546.6} & \textbf{5692.5} & \textbf{7248.8} & \textbf{5269.6} & \textbf{6855.1} & \texttt{124\_n2s27}$^\dagger$ & \textbf{19.0} & 32.9 & 21.4 & 20.0 & 26.3 & 19.2 \\
\texttt{018\_n3s62} & \textbf{(0.1\%)} & (0.1\%) & (0.2\%) & (0.1\%) & (0.2\%) & (0.1\%) & \texttt{125\_n2s100}$^\dagger$ & 1113.2 & 1238.1 & \textbf{777.3} & \textbf{649.2} & \textbf{771.3} & \textbf{787.8} \\
\texttt{019\_n3s122} & ($\infty$\%) & \textbf{(98.5\%)} & \textbf{(97.8\%)} & \textbf{(53.7\%)} & \textbf{(54.3\%)} & \textbf{(53.7\%)} & \texttt{127\_n2s34}$^\dagger$ & 35.3 & 37.8 & \textbf{26.6} & 41.8 & \textbf{32.2} & \textbf{30.3} \\
\texttt{020\_n3s101}$^\dagger$ & 240.1 & 268.0 & \textbf{137.3} & 286.8 & \textbf{237.6} & \textbf{235.9} & \texttt{128\_n2s27}$^\dagger$ & 13.9 & \textbf{10.5} & 15.2 & 15.1 & 16.2 & \textbf{12.0} \\
\texttt{021\_n4s71}$^\dagger$ & 211.6 & 247.7 & \textbf{208.0} & 242.1 & 259.0 & 230.1 & \texttt{129\_n2s65}$^\dagger$ & 786.4 & 1806.7 & 920.1 & \textbf{705.2} & 939.5 & \textbf{783.7} \\
\texttt{022\_n4s83}$^\dagger$ & 1365.4 & \textbf{1065.4} & \textbf{879.1} & \textbf{984.0} & \textbf{1124.7} & \textbf{977.2} & \texttt{130\_n2s85} & (1.5\%) & (80.1\%) & \textbf{(0.5\%)} & \textbf{(0.6\%)} & (1.8\%) & \textbf{(0.8\%)} \\
\texttt{023\_n4s97} & (0.3\%) & (0.7\%) & \textbf{(0.3\%)} & (0.3\%) & \textbf{(0.1\%)} & \textbf{(0.1\%)} & \texttt{131\_n3s96} & (0.1\%) & (0.1\%) & (0.1\%) & \textbf{(0.1\%)} & (0.1\%) & (0.1\%) \\
\texttt{025\_n4s93}$^\dagger$ & 1573.4 & 2580.3 & \textbf{1204.5} & 2210.2 & 1685.9 & 1915.1 & \texttt{132\_n3s56}$^\dagger$ & 837.0 & 3910.3 & \textbf{434.5} & \textbf{239.2} & \textbf{252.8} & 1239.2 \\
\texttt{026\_n4s111} & (13.4\%) & (14.1\%) & \textbf{(8.0\%)} & \textbf{(10.6\%)} & (17.4\%) & \textbf{(3.2\%)} & \texttt{133\_n3s76}$^\dagger$ & \textbf{35.4} & 75.4 & 42.6 & 43.3 & 57.7 & 50.7 \\
\texttt{027\_n4s122} & (22.4\%) & ($\infty$\%) & (88.7\%) & (46.0\%) & \textbf{(18.7\%)} & ($\infty$\%) & \texttt{134\_n3s62}$^\dagger$ & 295.3 & 306.4 & 366.8 & 316.0 & 296.6 & \textbf{265.8} \\
\texttt{028\_n4s168} & (91.5\%) & \textbf{(21.6\%)} & \textbf{(91.5\%)} & \textbf{(21.4\%)} & \textbf{(19.7\%)} & (91.6\%) & \texttt{135\_n3s55}$^\dagger$ & 313.7 & \textbf{257.9} & 342.1 & \textbf{308.1} & \textbf{267.2} & \textbf{233.6} \\
\texttt{030\_n4s177} & \textbf{(71.0\%)} & ($\infty$\%) & (72.0\%) & (72.0\%) & (72.0\%) & (72.0\%) & \texttt{136\_n3s63}$^\dagger$ & 1248.0 & 1343.3 & 1333.6 & \textbf{1034.4} & \textbf{1042.9} & \textbf{788.9} \\
\texttt{031\_n5s87} & (1.3\%) & \textbf{10488.0} & \textbf{(0.5\%)} & (8.6\%) & \textbf{9117.6} & \textbf{10642.5} & \texttt{137\_n3s97}$^\dagger$ & 6363.7 & \textbf{6014.9} & \textbf{5440.0} & 6591.0 & 7409.7 & \textbf{6213.3} \\
\texttt{032\_n5s88} & 8564.3 & (0.1\%) & (0.1\%) & \textbf{7273.1} & \textbf{6418.6} & (0.1\%) & \texttt{139\_n3s144} & (34.8\%) & \textbf{(17.1\%)} & (49.3\%) & (56.9\%) & \textbf{(9.5\%)} & \textbf{(30.2\%)} \\
\texttt{033\_n5s104}$^\dagger$ & 558.5 & 721.0 & 658.0 & 604.7 & \textbf{501.3} & 632.0 & \texttt{141\_n4s76}$^\dagger$ & \textbf{7.0} & 7.1 & 7.2 & 8.8 & 9.2 & 7.0 \\
\texttt{037\_n5s119} & (3.2\%) & (3.2\%) & \textbf{(1.6\%)} & \textbf{(0.2\%)} & \textbf{(0.1\%)} & \textbf{(0.1\%)} & \texttt{142\_n4s111}$^\dagger$ & 1814.3 & 2486.6 & \textbf{1697.2} & \textbf{1240.8} & \textbf{1277.3} & \textbf{1642.4} \\
\texttt{039\_n5s93}$^\dagger$ & 2660.7 & 2985.9 & \textbf{2351.8} & \textbf{2488.8} & \textbf{2270.7} & \textbf{2652.2} & \texttt{146\_n4s53} & 3646.7 & (1.2\%) & \textbf{2000.7} & 4237.9 & 9675.6 & (0.5\%) \\
\texttt{040\_n5s118} & (8.2\%) & (13.3\%) & (9.7\%) & (8.8\%) & (9.8\%) & \textbf{(5.9\%)} & \texttt{151\_n5s78}$^\dagger$ & 1890.0 & \textbf{1695.1} & \textbf{1011.5} & \textbf{1043.8} & \textbf{1448.8} & \textbf{1502.8} \\
\texttt{041\_n2s29}$^\dagger$ & \textbf{11.7} & 15.3 & 19.3 & 16.1 & 15.3 & 18.0 & \texttt{153\_n5s123} & (11.0\%) & (14.6\%) & (12.1\%) & (26.7\%) & (33.6\%) & \textbf{(9.2\%)} \\
\texttt{042\_n2s57}$^\dagger$ & 2923.1 & 3410.1 & \textbf{2756.5} & \textbf{2060.1} & \textbf{2666.7} & \textbf{2237.8} & \texttt{155\_n5s162} & --- & --- & (98.3\%) & (98.3\%) & (98.3\%) & (98.3\%) \\
\texttt{043\_n2s16}$^\dagger$ & 1.9 & \textbf{1.9} & 1.9 & 1.9 & 2.1 & 2.0 & \texttt{160\_n5s139} & (9.7\%) & \textbf{(4.6\%)} & \textbf{(3.8\%)} & \textbf{(2.9\%)} & \textbf{(3.6\%)} & \textbf{(3.4\%)} \\
\texttt{044\_n2s38} & (0.2\%) & (0.2\%) & \textbf{(0.2\%)} & \textbf{9959.8} & \textbf{(0.2\%)} & \textbf{9656.2} & \texttt{161\_n2s28}$^\dagger$ & 21.5 & 25.4 & \textbf{20.4} & 22.6 & 24.1 & 27.9 \\
\texttt{045\_n2s107} & \textbf{(6.0\%)} & (8.4\%) & (6.0\%) & (7.8\%) & (7.0\%) & (9.5\%) & \texttt{162\_n2s56}$^\dagger$ & 858.8 & 1147.3 & \textbf{730.2} & 876.9 & 979.2 & 889.5 \\
\texttt{047\_n2s20}$^\dagger$ & 5.1 & \textbf{3.5} & 7.1 & 5.8 & \textbf{5.0} & 7.2 & \texttt{164\_n2s71}$^\dagger$ & 160.9 & 169.2 & \textbf{114.8} & \textbf{130.0} & \textbf{138.4} & \textbf{142.9} \\
\texttt{048\_n2s82} & (47.9\%) & \textbf{(45.8\%)} & \textbf{(4.1\%)} & \textbf{(42.6\%)} & \textbf{(47.5\%)} & \textbf{(10.3\%)} & \texttt{166\_n2s43}$^\dagger$ & \textbf{62.4} & 80.0 & 71.8 & 66.4 & 81.1 & 75.6 \\
\texttt{049\_n2s42}$^\dagger$ & 82.8 & \textbf{59.5} & \textbf{32.9} & \textbf{44.9} & \textbf{40.5} & \textbf{36.5} & \texttt{167\_n2s50}$^\dagger$ & 522.7 & \textbf{493.3} & \textbf{332.6} & \textbf{346.4} & 528.9 & \textbf{376.0} \\
\texttt{050\_n2s113} & \textbf{(97.9\%)} & (97.9\%) & (97.9\%) & (97.9\%) & (97.9\%) & (97.9\%) & \texttt{168\_n2s65}$^\dagger$ & 3554.6 & \textbf{3095.8} & 4098.3 & 5726.2 & 7746.7 & 5621.0 \\
\texttt{051\_n3s67}$^\dagger$ & \textbf{3568.2} & 3811.0 & 3872.1 & 4076.4 & 3691.6 & 4487.0 & \texttt{169\_n2s47}$^\dagger$ & 230.9 & 262.8 & 245.0 & \textbf{179.0} & \textbf{151.2} & \textbf{211.6} \\
\texttt{052\_n3s62}$^\dagger$ & 3068.9 & \textbf{2422.7} & \textbf{2206.0} & \textbf{2765.8} & \textbf{1937.4} & \textbf{2106.3} & \texttt{170\_n2s64} & (0.4\%) & (1.1\%) & \textbf{8721.0} & \textbf{9421.3} & \textbf{10716.7} & \textbf{8293.7} \\
\texttt{053\_n3s68}$^\dagger$ & 384.0 & 871.4 & \textbf{333.7} & 1364.4 & 497.5 & 555.2 & \texttt{171\_n3s74}$^\dagger$ & 4766.6 & 9192.2 & 5170.5 & 4797.5 & 6276.2 & \textbf{3826.1} \\
\texttt{054\_n3s124} & (23.1\%) & \textbf{(15.4\%)} & \textbf{(15.5\%)} & \textbf{(19.4\%)} & (50.2\%) & \textbf{(4.2\%)} & \texttt{173\_n3s32}$^\dagger$ & 11.0 & \textbf{10.9} & 12.9 & \textbf{10.4} & \textbf{10.9} & \textbf{9.1} \\
\texttt{055\_n3s87}$^\dagger$ & 3145.7 & 3937.1 & 3188.9 & 4113.0 & 3704.7 & \textbf{2976.5} & \texttt{174\_n3s94} & 8122.8 & 9405.9 & \textbf{7483.8} & (0.1\%) & \textbf{4583.5} & 9118.3 \\
\texttt{056\_n3s134} & (98.1\%) & \textbf{(98.1\%)} & \textbf{(98.0\%)} & (98.1\%) & (98.1\%) & \textbf{(98.0\%)} & \texttt{175\_n3s70} & \textbf{(3.1\%)} & (6.3\%) & (4.3\%) & (5.4\%) & (4.3\%) & (4.6\%) \\
\texttt{057\_n3s188} & --- & (99.5\%) & (99.5\%) & (99.5\%) & (99.5\%) & (99.5\%) & \texttt{177\_n3s155} & (30.2\%) & (85.2\%) & \textbf{(25.6\%)} & (85.4\%) & (85.4\%) & (85.2\%) \\
\texttt{059\_n3s74}$^\dagger$ & 8959.5 & 9412.8 & \textbf{7684.7} & 10216.9 & 9684.8 & \textbf{7801.5} & \texttt{179\_n3s100} & (4.4\%) & \textbf{(1.5\%)} & \textbf{(2.5\%)} & \textbf{(3.3\%)} & \textbf{(1.8\%)} & \textbf{(2.1\%)} \\
\texttt{060\_n3s43}$^\dagger$ & 77.7 & \textbf{69.7} & \textbf{45.9} & \textbf{44.5} & \textbf{38.3} & \textbf{29.8} & \texttt{180\_n3s42}$^\dagger$ & \textbf{17.1} & 19.1 & 21.0 & 25.6 & 28.9 & 28.7 \\
\texttt{061\_n4s60}$^\dagger$ & 64.7 & \textbf{48.2} & 67.0 & 72.8 & 88.2 & \textbf{48.7} & \texttt{181\_n4s77}$^\dagger$ & 2092.8 & \textbf{1304.1} & \textbf{1220.4} & \textbf{1536.9} & \textbf{2085.8} & \textbf{1474.3} \\
\texttt{062\_n4s101} & (0.5\%) & \textbf{(0.2\%)} & \textbf{9909.2} & \textbf{8593.3} & \textbf{7744.2} & \textbf{(0.1\%)} & \texttt{182\_n4s97} & (1.7\%) & (32.9\%) & \textbf{(0.9\%)} & \textbf{(1.5\%)} & \textbf{(0.3\%)} & \textbf{(0.1\%)} \\
\texttt{063\_n4s133} & (5.3\%) & \textbf{(0.5\%)} & \textbf{(0.5\%)} & \textbf{(1.0\%)} & \textbf{(1.3\%)} & \textbf{(0.1\%)} & \texttt{183\_n4s113} & (1.0\%) & (1.3\%) & \textbf{4369.1} & \textbf{9448.8} & \textbf{(0.8\%)} & \textbf{9842.2} \\
\texttt{064\_n4s129} & (22.1\%) & (100.0\%) & (100.0\%) & (100.0\%) & \textbf{(21.8\%)} & (100.0\%) & \texttt{184\_n4s129} & (65.4\%) & \textbf{(60.4\%)} & \textbf{(55.7\%)} & \textbf{(62.9\%)} & \textbf{(64.8\%)} & \textbf{(61.8\%)} \\
\texttt{066\_n4s126} & ($\infty$\%) & ($\infty$\%) & \textbf{(74.5\%)} & ($\infty$\%) & \textbf{(74.8\%)} & \textbf{(74.5\%)} & \texttt{191\_n5s94}$^\dagger$ & 57.6 & \textbf{45.8} & \textbf{30.7} & \textbf{56.2} & \textbf{34.3} & \textbf{34.0} \\
\texttt{071\_n5s113} & (4.7\%) & \textbf{(2.2\%)} & \textbf{(3.8\%)} & \textbf{(1.9\%)} & \textbf{(0.8\%)} & \textbf{(2.6\%)} & \texttt{192\_n5s150} & (99.0\%) & --- & \textbf{(29.5\%)} & \textbf{(29.4\%)} & \textbf{(29.4\%)} & \textbf{(29.5\%)} \\
\texttt{072\_n5s183} & (3.3\%) & \textbf{(1.9\%)} & \textbf{(1.8\%)} & \textbf{(0.2\%)} & --- & --- & \texttt{193\_n5s138} & ($\infty$\%) & ($\infty$\%) & \textbf{(89.1\%)} & \textbf{(88.7\%)} & \textbf{(89.3\%)} & \textbf{(89.1\%)} \\
\multicolumn{7}{@{}c@{\quad}|@{\quad}}{\scalebox{0.6}{$\vdots$}} & \texttt{194\_n5s108}$^\dagger$ & 6624.3 & \textbf{6288.6} & \textbf{4406.8} & \textbf{6269.1} & \textbf{5006.6} & \textbf{5079.5} \\
\bottomrule
 &  &  &  &  &  &  & Solved & 86 & 82 & 86 & 89 & 88 & 87 \\
 &  &  &  &  &  &  & Solved by all$^\dagger$ & 79 & 79 & 79 & 79 & 79 & 79 \\
 &  &  &  &  &  &  & Sum$^\dagger$ & 100941.2 & 112118.3 & 86715.2 & 96097.3 & 100536.4 & 92977.6 \\
 &  &  &  &  &  &  & Mean$^\dagger$ & 1277.7 & 1419.2 & 1097.7 & 1216.4 & 1272.6 & 1176.9 \\
 &  &  &  &  &  &  & Sh. Geom. Mean$^\dagger$ & 300.5 & 334.8 & 270.8 & 285.2 & 288.3 & 275.3 \\
\cmidrule(lr){8-14}
 &  &  &  &  &  &  & Wilcoxon $p$ &  & 0.9996 & 0.0130 & 0.0136 & 0.2080 & 0.0081 \\
 &  &  &  &  &  &  & Significance &  & n.s. & * & * & n.s. & ** \\
\cmidrule(lr){8-14}
\end{tabular}
}
\end{table}

\clearpage

\section{Detailed Convergence Quality}
\label{sec:appendix-pd-integrals}

\begin{table}[htbp!]
\centering
\caption{Primal-dual integral values (billions, $\int_0^T (\text{primal}(t) - \text{dual}(t)) dt$, cost$\cdot$seconds). Lower is better. Bold: better than Standard. Wilcoxon test (one-sided). Significance: * $p<0.05$, ** $p<0.01$, *** $p<0.001$. Summary statistics (bottom right) are computed over both columns.}
\label{tab:pd-integrals}
\scalebox{0.5}{
\begin{tabular}{@{}lr@{\,}r@{\,}r@{\,}r@{\,}r@{\,}r@{\quad}|@{\quad}lr@{\,}r@{\,}r@{\,}r@{\,}r@{\,}r@{}}
\toprule
 & \multicolumn{2}{c}{\textbf{Baseline}} & \multicolumn{4}{c}{\textbf{ML-Based}} & & \multicolumn{2}{c}{\textbf{Baseline}} & \multicolumn{4}{c}{\textbf{ML-Based}} \\
\cmidrule(lr){2-3} \cmidrule(lr{2em}){4-7} \cmidrule(lr){9-10} \cmidrule(lr){11-14}
Instance & \makecell{Standard} & \makecell{\texttt{Rand.}\\\texttt{0.5P}} & \makecell{\texttt{ML}\\\texttt{0.5T-5C}} & \makecell{\texttt{ML}\\\texttt{0.5P}} & \makecell{\texttt{ML}\\\texttt{0.7P}} & \makecell{\texttt{ML}\\\texttt{0.1T-0.6P}\\\texttt{10C}} & Instance & \makecell{Standard} & \makecell{\texttt{Rand.}\\\texttt{0.5P}} & \makecell{\texttt{ML}\\\texttt{0.5T-5C}} & \makecell{\texttt{ML}\\\texttt{0.5P}} & \makecell{\texttt{ML}\\\texttt{0.7P}} & \makecell{\texttt{ML}\\\texttt{0.1T-0.6P}\\\texttt{10C}} \\
\midrule
\texttt{abilene} & 50.0 & \textbf{40.0} & \textbf{50.0} & \textbf{50.0} & \textbf{50.0} & \textbf{40.0} & \multicolumn{7}{c@{}}{\scalebox{0.6}{$\vdots$}} \\
\texttt{atlanta} & 257 & 505 & 317 & 269 & 377 & 317 & \texttt{076\_n5s122} & 14561 & 23101 & \textbf{13110} & \textbf{12832} & \textbf{13272} & \textbf{12702} \\
\texttt{di-yuan} & 10.1 & 10.3 & 10.1 & \textbf{10.1} & 20.1 & 20.1 & \texttt{081\_n2s160} & 46350 & \textbf{42150} & \textbf{39640} & \textbf{43530} & \textbf{44010} & \textbf{39440} \\
\texttt{france} & 1620 & \textbf{1260} & \textbf{1050} & \textbf{1070} & \textbf{1120} & \textbf{1060} & \texttt{083\_n2s76} & 22340 & 94287 & 31170 & 38010 & 35550 & 39799 \\
\texttt{geant} & 2420 & 5210 & 3370 & 4280 & 3590 & 3700 & \texttt{084\_n2s38} & 50.0 & \textbf{40.0} & 50.0 & \textbf{40.0} & \textbf{40.0} & \textbf{40.0} \\
\texttt{newyork} & 61.2 & 61.6 & 61.7 & 61.4 & 61.3 & 71.6 & \texttt{085\_n2s83} & 150 & 150 & 160 & 160 & 160 & 160 \\
\texttt{nobel-eu} & 1151 & 1251 & 1191 & 1161 & 1161 & 1171 & \texttt{086\_n2s66} & 4000 & 11760 & 4770 & 5710 & 6500 & 4730 \\
\texttt{nobel-germany} & 40.0 & \textbf{40.0} & 40.0 & 40.0 & 40.0 & \textbf{40.0} & \texttt{087\_n2s86} & 73647 & 112043 & \textbf{70107} & \textbf{62328} & \textbf{67328} & \textbf{72647} \\
\texttt{pdh} & 10.3 & 20.3 & 10.4 & 10.3 & 20.2 & 10.3 & \texttt{088\_n2s55} & 360 & \textbf{230} & \textbf{250} & \textbf{270} & \textbf{290} & \textbf{280} \\
\texttt{polska} & 10.0 & 20.0 & 20.0 & 20.0 & 20.0 & 20.0 & \texttt{091\_n3s62} & 190 & 200 & 210 & 210 & 210 & 210 \\
\texttt{sun} & 110 & 150 & 200 & 230 & 230 & 220 & \texttt{092\_n3s26} & 20.0 & \textbf{10.00} & 20.0 & 20.0 & 20.0 & 20.0 \\
\texttt{001\_n2s26} & 20.0 & \textbf{20.0} & \textbf{20.0} & \textbf{20.0} & \textbf{20.0} & \textbf{20.0} & \texttt{093\_n3s90} & 7369 & \textbf{6119} & \textbf{3040} & \textbf{4889} & \textbf{6769} & \textbf{4929} \\
\texttt{003\_n2s76} & 1040 & 1290 & \textbf{670} & \textbf{660} & \textbf{660} & \textbf{670} & \texttt{094\_n3s69} & 390 & \textbf{360} & \textbf{280} & \textbf{290} & \textbf{280} & \textbf{280} \\
\texttt{004\_n2s33} & 90.0 & \textbf{50.0} & \textbf{40.0} & \textbf{40.0} & \textbf{40.0} & \textbf{40.0} & \texttt{095\_n3s55} & 40.0 & \textbf{40.0} & 50.0 & 50.0 & \textbf{40.0} & 50.0 \\
\texttt{006\_n2s51} & 340 & \textbf{300} & \textbf{270} & \textbf{280} & \textbf{290} & \textbf{270} & \texttt{098\_n3s94} & 4681 & \textbf{3843} & 4811 & \textbf{4541} & \textbf{4591} & \textbf{4561} \\
\texttt{007\_n2s43} & 2090 & 5260 & \textbf{1940} & 2270 & \textbf{1770} & 2680 & \texttt{099\_n3s102} & 118161 & \textbf{114103} & \textbf{67746} & \textbf{69756} & \textbf{71684} & \textbf{67687} \\
\texttt{008\_n2s156} & 15540 & \textbf{15100} & \textbf{15160} & \textbf{15130} & \textbf{15060} & \textbf{15200} & \texttt{104\_n4s152} & 4716 & 5517 & 5167 & 5128 & 5148 & 5158 \\
\texttt{009\_n2s35} & 30.0 & 30.0 & 30.0 & 30.0 & 30.0 & 30.0 & \texttt{112\_n5s78} & 330 & \textbf{250} & \textbf{270} & \textbf{310} & \textbf{330} & \textbf{300} \\
\texttt{010\_n2s42} & 130 & \textbf{120} & \textbf{120} & \textbf{120} & \textbf{120} & \textbf{120} & \texttt{113\_n5s131} & 26668 & 30337 & \textbf{16443} & \textbf{17083} & \textbf{17723} & \textbf{17062} \\
\texttt{011\_n3s64} & 380 & 460 & 410 & 430 & 450 & 420 & \texttt{115\_n5s150} & 2037 & \textbf{1786} & \textbf{1802} & \textbf{1798} & \textbf{1799} & \textbf{1790} \\
\texttt{012\_n3s78} & 1820 & \textbf{1650} & \textbf{1320} & \textbf{1370} & \textbf{1410} & \textbf{1400} & \texttt{121\_n2s27} & 80.0 & \textbf{70.0} & 80.0 & 80.0 & 80.0 & 80.0 \\
\texttt{013\_n3s87} & 610 & \textbf{600} & 620 & 610 & 620 & 610 & \texttt{122\_n2s54} & 460 & 620 & 550 & 560 & 540 & 560 \\
\texttt{014\_n3s108} & 38702 & \textbf{35792} & 41271 & 47371 & 65469 & 56281 & \texttt{123\_n2s90} & 540 & 1150 & 760 & 900 & 840 & 770 \\
\texttt{015\_n3s91} & 640 & 670 & \textbf{570} & \textbf{630} & 650 & \textbf{640} & \texttt{124\_n2s27} & 100.0 & 140 & 100.0 & 100.0 & 100.0 & \textbf{90.0} \\
\texttt{017\_n3s84} & 7380 & 10430 & \textbf{7200} & \textbf{6890} & 7820 & \textbf{6980} & \texttt{125\_n2s100} & 4940 & 7770 & \textbf{2690} & \textbf{2700} & \textbf{2580} & \textbf{2490} \\
\texttt{018\_n3s62} & 290 & \textbf{280} & 290 & 290 & 290 & 290 & \texttt{127\_n2s34} & 100.0 & \textbf{90.0} & \textbf{90.0} & 101 & 120 & 110 \\
\texttt{020\_n3s101} & 150 & 190 & 160 & 160 & 160 & 160 & \texttt{128\_n2s27} & 20.0 & \textbf{20.0} & \textbf{20.0} & 20.0 & \textbf{20.0} & \textbf{20.0} \\
\texttt{021\_n4s71} & 1650 & \textbf{1410} & \textbf{1500} & 1670 & 1880 & 1780 & \texttt{129\_n2s65} & 2628 & \textbf{1269} & \textbf{1559} & \textbf{1429} & \textbf{1509} & \textbf{1499} \\
\texttt{022\_n4s83} & 4040 & 5120 & \textbf{3400} & \textbf{3880} & \textbf{3910} & \textbf{3740} & \texttt{130\_n2s85} & 74930 & 110020 & \textbf{49200} & \textbf{72410} & 91880 & \textbf{52010} \\
\texttt{023\_n4s97} & 2280 & \textbf{2181} & \textbf{1580} & \textbf{1620} & \textbf{1600} & \textbf{1560} & \texttt{131\_n3s96} & 2260 & \textbf{1880} & \textbf{2150} & \textbf{2120} & \textbf{2230} & \textbf{2150} \\
\texttt{025\_n4s93} & 3820 & 4700 & 4940 & 5270 & 4230 & 5220 & \texttt{132\_n3s56} & 100 & \textbf{100} & 110 & 110 & 110 & 110 \\
\texttt{026\_n4s111} & 13522 & 14082 & \textbf{11082} & 14092 & \textbf{12512} & 13601 & \texttt{133\_n3s76} & 160 & 480 & 370 & 160 & 160 & 230 \\
\texttt{027\_n4s122} & 76867 & --- & 90699 & 96654 & 98861 & --- & \texttt{134\_n3s62} & 80.0 & \textbf{70.0} & 90.1 & 90.0 & 90.0 & 90.0 \\
\texttt{028\_n4s168} & 111379 & \textbf{103689} & 120349 & \textbf{111179} & 115909 & \textbf{107999} & \texttt{135\_n3s55} & 290 & \textbf{220} & \textbf{200} & \textbf{220} & \textbf{220} & \textbf{220} \\
\texttt{030\_n4s177} & 107090 & --- & \textbf{104530} & \textbf{104740} & \textbf{104460} & \textbf{104400} & \texttt{136\_n3s63} & 486 & \textbf{454} & \textbf{453} & \textbf{453} & \textbf{464} & \textbf{463} \\
\texttt{031\_n5s87} & 13372 & \textbf{9722} & \textbf{6992} & \textbf{8102} & \textbf{8232} & \textbf{7712} & \texttt{137\_n3s97} & 8521 & 19598 & \textbf{7650} & \textbf{7732} & \textbf{7900} & \textbf{7680} \\
\texttt{032\_n5s88} & 526 & \textbf{486} & 527 & 537 & 526 & 527 & \texttt{139\_n3s144} & 18190 & 28090 & 29790 & 31370 & 31520 & 30350 \\
\texttt{033\_n5s104} & 869 & \textbf{819} & \textbf{650} & \textbf{829} & \textbf{849} & \textbf{769} & \texttt{141\_n4s76} & 20.0 & 20.0 & 20.0 & 20.0 & 20.0 & 20.0 \\
\texttt{037\_n5s119} & 2406 & 2590 & \textbf{2395} & 2443 & 2533 & \textbf{2361} & \texttt{142\_n4s111} & 1530 & \textbf{1480} & \textbf{900} & \textbf{970} & \textbf{980} & \textbf{930} \\
\texttt{039\_n5s93} & 4369 & \textbf{3889} & \textbf{2969} & \textbf{2699} & \textbf{2869} & \textbf{3109} & \texttt{146\_n4s53} & 110 & \textbf{100} & \textbf{60.0} & \textbf{60.0} & \textbf{60.0} & \textbf{60.0} \\
\texttt{040\_n5s118} & 3255 & 3265 & 3404 & 3284 & 3275 & 3274 & \texttt{151\_n5s78} & 3410 & 3890 & 3630 & 4320 & 5089 & 5389 \\
\texttt{041\_n2s29} & 20.0 & \textbf{10.0} & 20.0 & 20.0 & 20.0 & 20.0 & \texttt{153\_n5s123} & 35990 & 46660 & 36640 & 44080 & \textbf{35630} & 40210 \\
\texttt{042\_n2s57} & 2100 & 2840 & \textbf{1630} & \textbf{1760} & 2200 & \textbf{1610} & \texttt{160\_n5s139} & 21880 & 22750 & \textbf{14100} & \textbf{14111} & \textbf{14401} & \textbf{13851} \\
\texttt{043\_n2s16} & 10.00 & 10.00 & 10.00 & 10.00 & 10.00 & 10.00 & \texttt{161\_n2s28} & 20.1 & \textbf{20.1} & 20.1 & \textbf{20.1} & 20.1 & 20.1 \\
\texttt{044\_n2s38} & 310 & \textbf{290} & \textbf{270} & \textbf{270} & \textbf{280} & \textbf{280} & \texttt{162\_n2s56} & 170 & 180 & \textbf{170} & \textbf{170} & \textbf{170} & \textbf{170} \\
\texttt{045\_n2s107} & 10330 & 12750 & \textbf{9560} & 11960 & 12540 & 13060 & \texttt{164\_n2s71} & 330 & \textbf{270} & \textbf{310} & \textbf{260} & \textbf{250} & \textbf{260} \\
\texttt{047\_n2s20} & 10.0 & \textbf{10.0} & 10.0 & 10.0 & \textbf{10.0} & \textbf{10.0} & \texttt{166\_n2s43} & 66.2 & 71.8 & 77.7 & 76.7 & \textbf{57.7} & 70.4 \\
\texttt{048\_n2s82} & 861 & 881 & 941 & 941 & 971 & 941 & \texttt{167\_n2s50} & 230 & 230 & \textbf{180} & \textbf{180} & \textbf{180} & \textbf{180} \\
\texttt{049\_n2s42} & 110 & \textbf{80.0} & \textbf{70.0} & \textbf{70.0} & \textbf{70.0} & \textbf{70.0} & \texttt{168\_n2s65} & 21386 & \textbf{12472} & \textbf{7336} & \textbf{9234} & \textbf{12273} & \textbf{8555} \\
\texttt{050\_n2s113} & 107997 & \textbf{55218} & \textbf{107997} & \textbf{107997} & \textbf{107997} & \textbf{107997} & \texttt{169\_n2s47} & 108 & \textbf{106} & \textbf{103} & 108 & 108 & \textbf{103} \\
\texttt{051\_n3s67} & 1431 & \textbf{1381} & \textbf{1171} & \textbf{1141} & \textbf{1131} & \textbf{1151} & \texttt{170\_n2s64} & 4160 & 4510 & \textbf{4070} & \textbf{4090} & \textbf{4100} & \textbf{4080} \\
\texttt{052\_n3s62} & 760 & \textbf{650} & \textbf{680} & \textbf{670} & \textbf{670} & \textbf{670} & \texttt{171\_n3s74} & 18840 & \textbf{9854} & \textbf{3360} & \textbf{3372} & \textbf{14971} & \textbf{3380} \\
\texttt{053\_n3s68} & 20.0 & 30.0 & 30.0 & 30.0 & 30.0 & 30.0 & \texttt{173\_n3s32} & 20.0 & \textbf{20.0} & 50.0 & 40.0 & 50.0 & 40.0 \\
\texttt{054\_n3s124} & 10790 & 19590 & \textbf{8900} & \textbf{7830} & \textbf{7980} & \textbf{7920} & \texttt{174\_n3s94} & 9580 & 11150 & 10690 & \textbf{8031} & \textbf{8890} & \textbf{7070} \\
\texttt{055\_n3s87} & 18060 & 29099 & \textbf{13600} & \textbf{13700} & \textbf{11720} & \textbf{11310} & \texttt{175\_n3s70} & 230 & \textbf{220} & 240 & 240 & 240 & 230 \\
\texttt{056\_n3s134} & 130329 & \textbf{112589} & \textbf{108679} & \textbf{112719} & \textbf{108029} & \textbf{107999} & \texttt{177\_n3s155} & 45000 & 52841 & 52200 & \textbf{36801} & \textbf{38121} & \textbf{33871} \\
\texttt{059\_n3s74} & 14260 & \textbf{11970} & 15500 & \textbf{14140} & \textbf{14130} & 14750 & \texttt{179\_n3s100} & 8942 & \textbf{7389} & \textbf{6299} & \textbf{6214} & \textbf{6264} & \textbf{6248} \\
\texttt{060\_n3s43} & 32.3 & \textbf{31.9} & \textbf{31.9} & 32.4 & 32.8 & 32.4 & \texttt{180\_n3s42} & 150 & \textbf{140} & 150 & 150 & 190 & 180 \\
\texttt{061\_n4s60} & 350 & 370 & 420 & 410 & \textbf{330} & \textbf{280} & \texttt{181\_n4s77} & 6520 & 8990 & \textbf{3200} & \textbf{3320} & \textbf{4940} & \textbf{4990} \\
\texttt{062\_n4s101} & 41658 & \textbf{27029} & \textbf{17409} & \textbf{32378} & \textbf{22779} & \textbf{18429} & \texttt{182\_n4s97} & 14734 & \textbf{12227} & \textbf{7152} & \textbf{8243} & \textbf{8643} & \textbf{8163} \\
\texttt{063\_n4s133} & 57388 & 58637 & \textbf{52848} & \textbf{45438} & \textbf{45868} & \textbf{44688} & \texttt{183\_n4s113} & 14230 & \textbf{12320} & \textbf{13420} & \textbf{11920} & \textbf{12480} & 15869 \\
\texttt{064\_n4s129} & 109379 & 114018 & 116088 & 112588 & \textbf{94289} & 110188 & \texttt{184\_n4s129} & 58872 & \textbf{48052} & \textbf{24433} & \textbf{30634} & \textbf{31915} & \textbf{29114} \\
\texttt{071\_n5s113} & 35099 & 47168 & \textbf{32079} & 49368 & \textbf{25529} & 45598 & \texttt{191\_n5s94} & 340 & \textbf{230} & \textbf{240} & \textbf{310} & \textbf{300} & \textbf{280} \\
\texttt{072\_n5s183} & 58790 & \textbf{51821} & \textbf{48922} & \textbf{49012} & --- & --- & \texttt{192\_n5s150} & 97242 & --- & \textbf{63757} & \textbf{63677} & \textbf{63707} & \textbf{63637} \\
\texttt{075\_n5s142} & 62244 & \textbf{50186} & \textbf{28537} & \textbf{40166} & \textbf{38066} & \textbf{35747} & \texttt{194\_n5s108} & 1225 & \textbf{1114} & \textbf{993} & \textbf{1003} & \textbf{1014} & \textbf{1004} \\
\multicolumn{7}{@{}c@{\quad}|@{\quad}}{\scalebox{0.6}{$\vdots$}} &  &  &  &  &  &  &  \\
\bottomrule
 &  &  &  &  &  &  & Sum & 1888056 & 1705096 & 1605588 & 1682330 & 1639049 & 1490068 \\
 &  &  &  &  &  &  & Mean & 14636 & 13533 & 12446 & 13041 & 12805 & 11733 \\
 &  &  &  &  &  &  & Sh. Geo. Mean & 1348.4 & 1273.8 & 1198.8 & 1235.0 & 1232.7 & 1151.7 \\
\cmidrule(lr){8-14}
 &  &  &  &  &  &  & Wilcoxon $p$ &  & 0.6092 & 0.0004 & 0.0006 & 0.0004 & 0.0008 \\
 &  &  &  &  &  &  & Significance &  & n.s. & *** & *** & *** & *** \\
\cmidrule(lr){8-14}
\end{tabular}
}
\end{table}

\end{document}